\documentclass[11pt,reqno]{amsart}

\usepackage[T1]{fontenc}
\usepackage{lmodern}

\usepackage[letterpaper,margin=1in]{geometry}
\usepackage{amsmath,amssymb,amsthm,mathtools}
\usepackage{bm}
\usepackage{enumitem}
\usepackage{microtype}
\usepackage{mathrsfs}
\usepackage{booktabs}
\usepackage{hyperref}
\hypersetup{colorlinks=true,linkcolor=blue,citecolor=blue,urlcolor=blue,pdftitle={Two-Scale Calculus for Riesz Fractional Gradients: Compactness, Realization, and Periodic Homogenization},pdfauthor={Hermann Douanla}}

\newcommand{\R}{\mathbb{R}}
\newcommand{\Z}{\mathbb{Z}}
\newcommand{\T}{\mathbb{T}}
\newcommand{\eps}{\varepsilon}
\newcommand{\supp}{\operatorname{supp}}
\newcommand{\divs}{\operatorname{div}^{s}}
\newcommand{\divy}{\operatorname{div}^{s}_{y}}
\newcommand{\grads}{\nabla^{s}}
\newcommand{\gradys}{\nabla^{s}_{y}}
\newcommand{\weakto}{\rightharpoonup}
\newcommand{\twoweak}{\stackrel{2}{\rightharpoonup}}
\newcommand{\twostrong}{\stackrel{2}{\longrightarrow}}
\newcommand{\Id}{\mathrm{Id}}
\newcommand{\dist}{\operatorname{dist}}

\newcommand{\Hperzero}{\dot H^s_{\#}(Y)}
\newcommand{\Ho}{H_0^s(\Omega)}
\newcommand{\N}{\mathbb{N}}
\newcommand{\Ahom}{A^{\rm hom}}
\newcommand{\Abar}{\overline A}
\newcommand{\twograds}{\boldsymbol{\nabla}^{\,s}_{2}}

\newtheorem{theorem}{Theorem}[section]
\newtheorem{maintheorem}{Theorem}

\newtheorem{proposition}[theorem]{Proposition}
\newtheorem{lemma}[theorem]{Lemma}
\newtheorem{corollary}[theorem]{Corollary}
\newtheorem{definition}[theorem]{Definition}

\theoremstyle{remark}
\newtheorem{remark}[theorem]{Remark}
\newtheorem{example}[theorem]{Example}

\title[Two-scale calculus for Riesz fractional gradients]{Two-Scale Calculus for Riesz Fractional Gradients: Compactness, Realization, and Periodic Homogenization}
\author[H. Douanla]{Hermann Douanla}
\date{}

\begin{document}

\thispagestyle{plain}
\markboth{H. DOUANLA}{TWO-SCALE CALCULUS FOR RIESZ FRACTIONAL GRADIENTS}

\vspace*{2.5\baselineskip}
\begin{center}
\parbox{0.86\textwidth}{\centering\bfseries
TWO-SCALE CALCULUS FOR RIESZ FRACTIONAL GRADIENTS: COMPACTNESS, REALIZATION, AND PERIODIC HOMOGENIZATION}

\vspace{1.6\baselineskip}
{\small HERMANN DOUANLA}\\[0.2\baselineskip]
Department of Mathematics, Faculty of Science, University of Yaound\'e I\\
P.O. Box 812, Yaound\'e, Cameroon\\
\texttt{hermann.douanla@facsciences-uy1.cm}
\end{center}

\vspace{0.8\baselineskip}
\begin{center}
\begin{minipage}{0.86\textwidth}
\small
\textbf{Abstract.} We develop a two-scale calculus for the Riesz fractional gradient
$\grads$, $0<s<1$, acting on functions that vanish outside a bounded
open set $\Omega\subset\R^N$. For bounded families in $H_0^s(\Omega)$
we characterize the interior two-scale limits of the fractional
gradients and prove, by a Fourier-multiplier estimate, that every limit
of this form is generated by such a family; a periodic fractional
Helmholtz decomposition identifies the microscopic potentials, and
$\eps^s$ is the natural amplitude of the oscillations. For periodic
fractional divergence-form equations with bounded measurable, uniformly
elliptic, not necessarily symmetric coefficients, we identify the two-scale limit of the
fractional gradients and the homogenized problem, whose coefficient is
the classical homogenized matrix in $\Omega$ and the cell average
outside $\Omega$, and we construct correctors in $H_0^s(\Omega)$ that
give strong convergence in $H^s(\R^N)$.
The identification extends to parabolic problems with dual-space
forcing.


\end{minipage}
\end{center}

\vspace{2\baselineskip}

\section{Introduction}

\subsection{Problem and motivation}

For $0<s<1$, the Riesz fractional gradient $\grads$ is the Fourier multiplier
\[
 \widehat{\grads u}(\xi)=i\xi|\xi|^{s-1}\widehat u(\xi)\qquad(\xi\in\R^N).
\]
It is singled out by its invariance properties: by a theorem of \v{S}ilhav\'y \cite[Theorem~2.2]{Silhavy2020}, every linear operator $G$ on $C_c^\infty(\R^N)$ with vector values that is invariant under translations and rotations, homogeneous of degree $s$ under dilations, and continuous in the sense that $f\mapsto Gf(x)$ is a distribution for each $x\in\R^N$, is a constant multiple of $\grads$. The Riesz fractional gradient has become a central object of the fractional calculus of variations \cite{ComiStefani2019,KreisbeckSchoenberger2022,ShiehSpector2015}. This paper develops a two-scale calculus for $\grads$; its main application is the complementary-value problem
\begin{equation}\label{eq:intromodel}
 -\divs\!\left(A\!\left(\frac{x}{\eps}\right)\grads u_\eps\right)=f
 \quad\text{in }\Omega,
 \qquad
 u_\eps=0\quad\text{in }\Omega^c,
\end{equation}
with $A$ periodic and uniformly elliptic; the prescription $u_\eps=0$ in $\Omega^c$ is the complementary-value condition. Since $\grads u$ is not supported in $\Omega$ even when $u$ is, the limit combines two mechanisms: in $\Omega$ the coefficient homogenizes to the classical matrix $\Ahom$, whereas in $\Omega^c$ it averages to $\Abar=\int_YA$, where $Y=(0,1)^N$ is the reference cell. Three questions are invisible at the level of the effective matrix: which microscopic limits can oscillating Riesz gradients have, is every such limit generated by a family satisfying the complementary-value condition, and at which amplitude does the microscopic potential enter? We answer them by a fractional two-scale calculus and an exact compactness--realization theory, and we use them to identify the two-scale limit of the fractional gradients, the homogenized problem, and admissible correctors, that is, correctors in $H_0^s(\Omega)$.

\subsection{Literature review}

Periodic homogenization and correctors are classical \cite{BensoussanLionsPapanicolaou1978,CioranescuDonato1999,JikovKozlovOleinik1994,Tartar2009}. Two-scale convergence originates with Nguetseng \cite{Nguetseng1989} and Allaire \cite{Allaire1992}. It has been extended to periodic hypersurfaces \cite{AllaireDamlamianHornung1996,NeussRadu1996}, reiterated homogenization \cite{AllaireBriane1996}, periodic measures and double porosity \cite{Zhikov2000}, stochastic settings \cite{BourgeatMikelicWright1994}, generalized Besicovitch spaces \cite{CasadoDiazGayte2002,SangoSvanstedtWoukeng2011}, homogenization algebras \cite{Nguetseng2025}, measures and $BV$ \cite{Amar1998}, and Orlicz--Sobolev spaces \cite{TachagoNnang2012}.

Visintin developed a two-scale calculus for differentiable functions \cite{Visintin2006} and a characterization and converse realization for classical first-order operators \cite[Theorem~2.2(ii)]{VisintinFirstOrder2007}; he also pointed to fractional Sobolev two-scale convergence as a possible direction \cite[Section~5]{Visintin2006}. Giorgio, Happ and Sch\"onberger \cite{GiorgioHappSchoenberger2026} introduced a two-scale convergence adapted to nonlocal micromagnetic exchange energies, in which the interaction range and the heterogeneity scale vanish together, and obtained a nonlocal analogue of the decomposition of two-scale limits of gradients for difference quotients \cite[Theorem~4.3]{GiorgioHappSchoenberger2026}; Braides, Brusca and Donati \cite{BraidesBruscaDonati2024} recover classical elliptic homogenization from Gagliardo seminorms with oscillating coefficients as $s\to1$. Other two-scale nonlocal models appear in \cite{DuLiptonMengesha2016,RamirezTorresPentaGrillo2022}. Here the order $s$ is fixed and the interactions have infinite range, as for the Riesz gradient itself.

On the operator side, up to a Fourier normalization, $\grads$ is the gradient $D^s$ of Shieh--Spector \cite{ShiehSpector2015,ShiehSpector2018}; its distributional and regularity theory is developed in \cite{ComiStefani2019,SchikorraShiehSpector2015}, the localization of nonlocal gradients with integrable kernels is studied in \cite{MengeshaSpector2015}, and transfer operators between classical, fractional and finite-horizon gradients are used by Cueto, Kreisbeck and Sch\"onberger \cite{CuetoKreisbeckSchoenberger2023} for truncated interactions.

For operators in fractional divergence form, Caponi, Carbotti and Maione \cite[Theorems~3.1 and 5.1]{CaponiCarbottiMaione2025} prove $H$-compactness with coefficients prescribed outside $\Omega$, show that the nonlocal $H$-limit coincides in $\Omega$ with the local one and relate it to $\Gamma$-convergence for symmetric coefficients. Buchinger et al.\ \cite[Theorem~1.3]{BuchingerEtAl2026}, under the assumption $|\partial\Omega|=0$, characterize $H^s$-convergence for uniformly elliptic coefficients by local $H$-convergence in $\Omega$ together with weak-* convergence of the coefficients on $\R^N\setminus\overline\Omega$, and compare $H^s$-convergence with the nonlocal $H$-convergence of \cite{Waurick2018}, an abstract framework for operator-valued coefficients. The nonlinear monotone case, and its relation to $\Gamma$-convergence, is treated in \cite{BruscaEtAl2026}. The broader nonlocal vector calculus underlying the fractional divergence and Helmholtz structure is discussed in \cite{DEliaGulianMengeshaScott2022,DEliaGulianOlsonKarniadakis2021,Silhavy2020}.

The periodic effective law of Theorem~\ref{thm:mainB} below agrees with this operator-level characterization (Remark~\ref{rem:effectiveknown}); here this law is obtained from the microscopic structure, which is what gives access to the two-scale limit of the gradients and to admissible correctors. Compactness can be obtained either directly, with the rate of Lemma~\ref{lem:oscdiv}, or by transfer, since the Riesz potential $I_{1-s}$, with multiplier $|\xi|^{s-1}$, maps $H_0^s(\Omega)$ into $H^1_{\rm loc}(\R^N)$ with $\nabla I_{1-s}=\grads$ (Remark~\ref{rem:transferproof}). Realization, the converse of compactness, does not transfer: a classical recovery family for $I_{1-s}u_0$ is in general not the Riesz potential of a function vanishing outside $\Omega$, and the fractional gradient of a compactly supported oscillation is not localized, so that its exterior part must be controlled in $L^2(\R^N)$ to keep the family bounded in $H^s(\R^N)$.

\subsection{Main results}

The new results of the paper are the exact characterization of the interior two-scale limits of Riesz gradients under the complementary-value condition, with converse realization at amplitude $\eps^s$, and the identification of the two-scale limit of the fractional gradients through fractional cell potentials. Without symmetry assumptions, we also prove separate convergence of the interior and exterior energies and construct strong admissible correctors; these yield strong convergence of $\grads u_\eps$ on all of $\Omega^c$, whereas general compactness results give it away from $\partial\Omega$ \cite{BruscaEtAl2026,KreisbeckSchoenberger2022}. Finally, we identify the time-dependent two-scale limit of the gradients for parabolic problems with dual-space forcing.

Theorem~\ref{thm:mainA} is the foundational result of the paper; Theorems~\ref{thm:mainB} and~\ref{thm:mainC} are its applications to elliptic and parabolic homogenization. Here $H_0^s(\Omega)$ is the closure of $C_c^\infty(\Omega)$ in $H^s(\R^N)$; its elements vanish a.e.\ outside $\Omega$ (Subsection~\ref{sec:fourier}). The space $\dot H^s_\#(Y)$ consists of the zero-mean $Y$-periodic $H^s$ functions (Subsection~\ref{sec:periodicspaces}), and the arrows $\twoweak$ and $\twostrong$ denote weak and strong two-scale convergence (Subsection~\ref{sec:twoscale}). The analysis takes place in the two-scale state space
\[
 \mathbb H_{2,0}^s(\Omega\times Y):=H_0^s(\Omega)\times L^2\bigl(\Omega;\dot H_\#^s(Y)\bigr),
 \qquad
 \twograds U:=\grads u_0+\gradys u_1\quad(U=(u_0,u_1)).
\]
A family $(u_\eps)_{\eps>0}$ \emph{realizes} $U=(u_0,u_1)$ if it is bounded in $H_0^s(\Omega)$, $u_\eps\to u_0$ in $L^2(\R^N)$, and $\grads u_\eps|_\Omega\twostrong\twograds U$ in $L^2(\Omega\times Y;\R^N)$.

\begin{maintheorem}\label{thm:mainA}
Let $0<s<1$ and let $\Omega\subset\R^N$ be bounded and open.
\begin{enumerate}[label=(\alph*)]
\item If $(u_\eps)_{\eps>0}$ is uniformly bounded in $H_0^s(\Omega)$, then for every null sequence $\eps_n\downarrow0$ there exist a subsequence, not relabeled, and a unique $U=(u_0,u_1)\in\mathbb H_{2,0}^s(\Omega\times Y)$ such that
\[
 u_{\eps_n}\to u_0\quad\text{in }L^2(\R^N),
 \qquad
 \grads u_{\eps_n}\twoweak\grads u_0+\gradys u_1
 \quad\text{in }L^2(\Omega\times Y;\R^N).
\]
\item Conversely, every $U\in\mathbb H_{2,0}^s(\Omega\times Y)$ is realized by a family $(u_\eps)_{\eps>0}\subset H_0^s(\Omega)$.
\item If $U=(u_0,u_1)$ with $u_1\in C_c^\infty(\Omega;C_\#^\infty(Y))$ of zero $Y$-mean, the family $u_\eps=u_0+\eps^su_1(x,x/\eps)$ realizes $U$, and
\[
 \bigl\|\grads u_\eps-\grads u_0
 -(\gradys u_1)(x,x/\eps)\bigr\|_{L^2(\R^N)}
 \le C(u_1,s)\eps^s.
\]
Within this fixed-profile ansatz, $s$ is the only exponent producing a bounded family with a nontrivial microscopic gradient (Proposition~\ref{prop:amplitudedichotomy}).
\end{enumerate}
\end{maintheorem}

Since realization does not transfer, Theorem~\ref{thm:mainA}(b)--(c) rests on a direct construction: a Fourier-multiplier estimate (Lemma~\ref{lem:amplitude}) shows that $\eps^su_1(x,x/\eps)$ has fractional gradient $(\gradys u_1)(x,x/\eps)+O(\eps^s)$ in $L^2(\R^N)$, which identifies $\eps^s$ as the natural amplitude, and a diagonal argument extends the realization to every $u_1\in L^2(\Omega;\dot H^s_\#(Y))$. Our direct proof of compactness uses a cancellation estimate, with rate $\eps^{1-s}$, for oscillating fields that are fractional-divergence-free in $y$ (Lemma~\ref{lem:oscdiv}), a mechanism that differs from the case $s=1$ (Remark~\ref{rem:oscdivs1}).

For the homogenization problems, the bridge $\mathcal J_s=(-\Delta_y)^{(1-s)/2}$, an isometry from $\dot H^1_\#(Y)$ onto $\dot H^s_\#(Y)$ with $\gradys\mathcal J_s=\nabla_y$ (Theorem~\ref{thm:bridge}), shows that the fractional and classical microscopic gradient spaces coincide while the potentials differ.

\begin{maintheorem}\label{thm:mainB}
Let $0<s<1$. Assume that $A\in L_\#^\infty(Y;\R^{N\times N})$ is uniformly elliptic, let $\Omega\subset\R^N$ be bounded and open with $|\partial\Omega|=0$, let $f\in(H_0^s(\Omega))'$, and let $u_\eps$ solve \eqref{eq:intromodel} in the variational sense. Then
\[
 u_\eps\rightharpoonup u_0\quad\text{in }H_0^s(\Omega),
 \qquad
 u_\eps\to u_0\quad\text{in }L^2(\R^N),
\]
where $u_0\in H_0^s(\Omega)$ is the unique solution of
\begin{equation}\label{eq:mainChom}
 \int_\Omega \Ahom\grads u_0\cdot\grads v\,dx
 +\int_{\Omega^c}\Abar\grads u_0\cdot\grads v\,dx
 =\langle f,v\rangle
 \qquad\forall v\in H_0^s(\Omega),
\end{equation}
with $\Abar=\int_YA(y)\,dy$ and $\Ahom$ the classical periodic homogenized matrix.

Moreover, the two-scale limit of the fractional gradients is identified:
\[
 \grads u_\eps\twoweak
 [\Id+Q(y)]\grads u_0(x)
 \quad\text{in }L^2(\Omega\times Y;\R^N),
\]
where $Q$ is the matrix with columns $\gradys\chi_1^s,\dots,\gradys\chi_N^s$, generated by the fractional cell potentials $\chi_j^s=\mathcal J_s\chi_j$, with $\chi_j\in\dot H^1_\#(Y)$ the classical cell correctors.

The energies converge separately in $\Omega$ and in $\Omega^c$:
\begin{align*}
 \int_\Omega A(x/\eps)\grads u_\eps\cdot\grads u_\eps\,dx&\to\int_\Omega\Ahom\grads u_0\cdot\grads u_0\,dx,\\
 \int_{\Omega^c}A(x/\eps)\grads u_\eps\cdot\grads u_\eps\,dx&\to\int_{\Omega^c}\Abar\grads u_0\cdot\grads u_0\,dx,
\end{align*}
and, with $A^{\rm sym}=\frac12(A+A^T)$ and $G(x,y)=[\Id+Q(y)]\grads u_0(x)$,
\[
 (A^{\rm sym})^{1/2}(x/\eps)\grads u_\eps\twostrong
 (A^{\rm sym})^{1/2}(y)G(x,y)
 \quad\text{in }L^2(\Omega\times Y;\R^N).
\]

There exist correctors $\mathcal C_\eps\in H_0^s(\Omega)$ such that
\[
 u_\eps-u_0-\mathcal C_\eps\to0\quad\text{in }H^s(\R^N).
\]
Moreover, $\grads\mathcal C_\eps\to0$ in $L^2(\Omega^c;\R^N)$; in particular, $\grads u_\eps\to\grads u_0$ strongly in $L^2(\Omega^c;\R^N)$.
\end{maintheorem}

\begin{maintheorem}\label{thm:mainC}
Let $0<s<1$, $T>0$, and let $\Omega\subset\R^N$ be bounded and open
with $|\partial\Omega|=0$. Assume that
$A\in L^\infty_\#(Y;\R^{N\times N})$ is uniformly elliptic, and that
\[
 f_\eps\to f\quad\text{in }L^2(0,T;(H_0^s(\Omega))'),
 \qquad
 u_\eps^{\rm in}\to u^{\rm in}\quad\text{in }L^2(\Omega).
\]
Let $u_\eps\in L^2(0,T;H_0^s(\Omega))\cap H^1(0,T;(H_0^s(\Omega))')$ be the variational solution of
\[
\begin{cases}
 \partial_t u_\eps-\divs\!\bigl(A(x/\eps)\grads u_\eps\bigr)=f_\eps
   &\text{in }(0,T)\times\Omega,\\
 u_\eps=0&\text{in }(0,T)\times\Omega^c,\\
 u_\eps(0,\cdot)=u_\eps^{\rm in}&\text{in }\Omega.
\end{cases}
\]
Define $B_{\rm eff}=\Ahom$ on $\Omega$ and $B_{\rm eff}=\Abar$
on $\Omega^c$, with $\Ahom$ and $\Abar$ as in Theorem~\ref{thm:mainB}. Then $u_\eps\to u_0$ strongly in
$L^2((0,T)\times\Omega)$, where $u_0$ is the unique variational
solution of
\[
\begin{cases}
 \partial_t u_0-\divs\!\bigl(B_{\rm eff}(x)\grads u_0\bigr)=f
   &\text{in }(0,T)\times\Omega,\\
 u_0=0&\text{in }(0,T)\times\Omega^c,\\
 u_0(0,\cdot)=u^{\rm in}&\text{in }\Omega.
\end{cases}
\]
Moreover,
\[
 \grads u_\eps\twoweak[\Id+Q(y)]\grads u_0(t,x)
 \quad\text{in }L^2((0,T)\times\Omega\times Y;\R^N),
\]
with $Q$ as in Theorem~\ref{thm:mainB}.
\end{maintheorem}

Both equations are understood in the variational sense of Section~\ref{sec:parabolicsub}, and the two-scale convergence in the space--time sense defined there.

The elliptic framework also applies to oscillatory flux forcing. In Section~\ref{sec:oscforcing}, an additional periodic cell potential identifies the effective forcing, which may be nonzero even when the oscillating forcing field has zero cell mean; as in the classical case, it is independent of $s$.

The two-scale calculus developed here is not tied to the model problems treated in this paper. It provides the compactness, realization and corrector tools needed for further homogenization problems for Riesz fractional gradients, such as high-contrast and double-porosity media, reiterated and stochastic homogenization, monotone nonlinear operators, fractional elasticity, and quantitative convergence rates.

\subsection{Organization of the paper}
Section~\ref{sec:prelim} collects the functional framework and Section~\ref{sec:helmholtz} the periodic fractional calculus. Theorem~\ref{thm:mainA} is proved in Section~\ref{sec:oscdiv}: part~(a) is Theorem~\ref{thm:maincompactness}, part~(b) is Theorem~\ref{thm:characterization}, and part~(c) is Theorem~\ref{thm:realization} with Proposition~\ref{prop:amplitudedichotomy}. Theorem~\ref{thm:mainB} is proved in Section~\ref{sec:homogenization}: it follows from Proposition~\ref{prop:classicalcell}, Theorems~\ref{thm:homogenization}, \ref{thm:strongenergy} and~\ref{thm:generalcorrector}, and Corollary~\ref{cor:exteriorstrong}. Section~\ref{sec:oscforcing} treats oscillatory flux forcing (Proposition~\ref{prop:oscforcing}), and Theorem~\ref{thm:mainC} is Theorem~\ref{thm:parabolic} (Section~\ref{sec:parabolicsub}).

\section{Preliminaries and notation}\label{sec:prelim}

Throughout the paper, $0<s<1$, $Y=(0,1)^N$, $\T^N=\R^N/\Z^N$, and $\Omega\subset\R^N$ is a bounded open set. From Section~\ref{sec:homogenization} on we also assume $|\partial\Omega|=0$; this is the only property of the boundary used in the homogenization and corrector results. Sections~\ref{sec:helmholtz} and \ref{sec:oscdiv} hold without it. We write $\mathbb B_t$ for the open ball of radius $t$ centered at $0$. The letter $C$ denotes a constant whose value may change from line to line and which does not depend on $\eps$.

\subsection{Fourier convention and Riesz fractional operators}\label{sec:fourier}

We use the unitary Fourier transform
\[
 \widehat u(\xi)=(2\pi)^{-N/2}\int_{\R^N}e^{-ix\cdot\xi}u(x)\,dx.
\]
For $u\in\mathscr S(\R^N)$, define
\[
 \widehat{\grads u}(\xi)=i\xi|\xi|^{s-1}\widehat u(\xi),\qquad\xi\in\R^N\setminus\{0\}.
\]
For vector fields $F\in\mathscr S(\R^N;\R^N)$,
\begin{equation}\label{eq:defdivs}
 \widehat{\divs F}(\xi)=i\xi\cdot\widehat F(\xi)|\xi|^{s-1},\qquad\xi\in\R^N\setminus\{0\}.
\end{equation}
For definiteness, let $\mathcal R^+$ denote the vector Riesz transform with multiplier
\[
 \widehat{\mathcal R^+ f}(\xi)=i\frac{\xi}{|\xi|}\widehat f(\xi).
\]
Then, with the Fourier convention above,
\[
 \grads=\mathcal R^+(-\Delta)^{s/2}=\nabla I_{1-s},
 \qquad
 \widehat{I_{1-s}f}(\xi)=|\xi|^{s-1}\widehat f(\xi).
\]
If instead the standard harmonic-analysis convention $\widehat{Rf}(\xi)=-i\xi|\xi|^{-1}\widehat f(\xi)$ is used, then $\grads=-R(-\Delta)^{s/2}$. Thus our operator is the Fourier-normalized form of the fractional gradient $D^s$ of Shieh--Spector; only the sign and the normalizing constant change with the normalization of the Riesz transform and potential. With these conventions,
\[
 -\divs\grads u=(-\Delta)^s u.
\]

The Bessel-potential Sobolev space is
\[
 H^s(\R^N)=\left\{u\in L^2(\R^N):\int_{\R^N}(1+|\xi|^2)^s|\widehat u(\xi)|^2\,d\xi<\infty\right\}.
\]
Its dual is denoted $H^{-s}(\R^N)$. For $u\in H^s(\R^N)$ and $F\in L^2(\R^N;\R^N)$,
\[
 \int_{\R^N}\grads u\cdot F\,dx
 =-\langle u,\divs F\rangle,
\]
where $\divs F\in H^{-s}(\R^N)$ is defined by \eqref{eq:defdivs} in the sense of distributions. We also use the homogeneous seminorm
\[
 \|u\|_{\dot H^s(\R^N)}^2:=\int_{\R^N}|\xi|^{2s}|\widehat u(\xi)|^2\,d\xi=\|\grads u\|_{L^2(\R^N)}^2 .
\]
Unless indicated otherwise, $\|\cdot\|_2$ denotes the norm of $L^2(\R^N)$. The space $C_\#(Y)$ consists of continuous $Y$-periodic functions on $\R^N$, with the supremum norm, and $C^\infty_\#(Y)$ of smooth ones.

For a bounded open set $\Omega$, we use the complementary-value fractional Dirichlet space
\begin{equation}\label{eq:H0s}
 \Ho:=\overline{C_c^\infty(\Omega)}^{\,H^s(\R^N)}.
\end{equation}
We write $\Omega^c:=\R^N\setminus\Omega$. Every element of $\Ho$ vanishes a.e.\ on $\Omega^c$. In the notation of the boundary-integral literature, $\Ho$ is the space $\widetilde H^s(\Omega)$; it should not be confused with the closure of $C_c^\infty(\Omega)$ in $H^s(\Omega)$. On $\Ho$ the norm $\|\grads u\|_{L^2(\R^N)}$ is equivalent to the full $H^s$ norm. Indeed, with the unitary Fourier transform and $u=0$ a.e.\ on $\Omega^c$, Cauchy--Schwarz yields
$|\widehat u(\xi)|^2\le (2\pi)^{-N}|\Omega|\|u\|_{L^2(\R^N)}^2$.
For any $R>0$, Plancherel and the definition of $\grads$ therefore give
\[
 \|u\|_2^2 \le (2\pi)^{-N}|\mathbb B_1|\,|\Omega|R^N\|u\|_2^2
       +R^{-2s}\|\grads u\|_2^2.
\]
Choose $R$ so that $(2\pi)^{-N}|\mathbb B_1|\,|\Omega|R^N\le\tfrac12$ and absorb the first term. Hence
\[
 \|u\|_{L^2(\R^N)}\le C(N,s,|\Omega|)\|\grads u\|_{L^2(\R^N)}
 \qquad (u\in\Ho).
\]
No boundary regularity is needed for this inequality. For the general fractional Poincar\'e framework, see \cite{KreisbeckSchoenberger2022,ShiehSpector2015}.

\subsection{Periodic fractional spaces}\label{sec:periodicspaces}

For $k\in\Z^N$ set $q_k=2\pi k$; throughout, $k\ne0$ means $k\in\Z^N\setminus\{0\}$. A $Y$-periodic function or distribution is written $v(y)=\sum_{k\in\Z^N}\widehat v(k)e^{iq_k\cdot y}$. For $1\le r\le\infty$, $L^r_\#(Y;\R^d)$ denotes the $Y$-periodic functions on $\R^N$ whose restriction to $Y$ lies in $L^r(Y;\R^d)$, and $L^2_{\#,0}(Y;\R^d)$ the subspace of zero-mean functions of $L^2_\#(Y;\R^d)$. For $t\in\R$, let $H^t_\#(Y)$ be the space of periodic distributions with $\sum_{k\in\Z^N}(1+|q_k|^2)^t|\widehat v(k)|^2<\infty$. Define
\[
 \Hperzero
 =\left\{v\in H^s_\#(Y):\int_Yv\,dy=0\right\},
\]
with norm
\[
 \|v\|_{\Hperzero}^2
 :=\sum_{k\in\Z^N\setminus\{0\}}|q_k|^{2s}|\widehat v(k)|^2.
\]
The periodic fractional gradient is
\[
 \widehat{\gradys v}(k)
 =iq_k|q_k|^{s-1}\widehat v(k),\qquad k\ne0,
\]
and the periodic fractional divergence is
\[
 \widehat{\divy F}(k)
 =iq_k\cdot\widehat F(k)|q_k|^{s-1},\qquad k\ne0.
\]
The space $\dot H^1_\#(Y)$ is defined in the same way, with the weight $|q_k|^{2}$ ($k\ne0$). By Parseval,
\begin{equation}\label{eq:isometry}
 \|\gradys v\|_{L^2(Y)}=\|v\|_{\Hperzero}.
\end{equation}

\subsection{Two-scale convergence}\label{sec:twoscale}

The parameter $\eps>0$ tends to zero. A notation such as $(v_\eps)_{\eps>0}\subset X$ denotes the corresponding $\eps$-indexed family; as is customary in homogenization, we occasionally refer to such families as sequences. Whenever compactness requires an actual subsequence, this means extraction along a null sequence $\eps_n\downarrow0$. A family, or a corrector, is called \emph{admissible} if it belongs to $\Ho$, that is, if it satisfies the complementary-value condition.

We use repeatedly the periodic averaging principle: for $w\in L^1_\#(Y)$, $w(x/\eps)\to\int_Yw\,dy$ weakly in $L^1_{\rm loc}(\R^N)$, and weakly-* in $L^\infty(\R^N)$ if $w\in L^\infty_\#(Y)$. We recall the periodic two-scale convergence of Nguetseng and Allaire \cite{Allaire1992,LukkassenNguetsengWall2002,Nguetseng1989}. For $d\ge1$ we use test functions in $L^2\bigl(\Omega;C_\#(Y;\R^d)\bigr)$, equipped with the Bochner norm induced by the supremum norm in $y$; for such $W$ we write $W^\eps(x):=W(x,x/\eps)$, $x\in\Omega$.

\begin{definition}
A family $(v_\eps)_{\eps>0}\subset L^2(\Omega;\R^d)$ converges weakly two-scale to $v_0\in L^2(\Omega\times Y;\R^d)$, written $v_\eps\twoweak v_0$, if
\[
 \lim_{\eps\to0}\int_\Omega v_\eps\cdot\varphi^\eps\,dx
 =\int_{\Omega\times Y}v_0\cdot\varphi\,dxdy
 \qquad\forall\varphi\in L^2(\Omega;C_\#(Y;\R^d)).
\]
It converges strongly two-scale, written $v_\eps\twostrong v_0$, if moreover $\|v_\eps\|_{L^2(\Omega)}\to\|v_0\|_{L^2(\Omega\times Y)}$.
\end{definition}

For uniformly bounded families it suffices to test with smooth $\varphi\in C_c^\infty(\Omega;C_\#^\infty(Y;\R^d))$. Every family bounded in $L^2(\Omega;\R^d)$ has, along every null sequence, a weakly two-scale convergent subsequence, and $v_\eps\twoweak v_0$ implies $v_\eps\weakto\int_Yv_0(\cdot,y)\,dy$ in $L^2(\Omega;\R^d)$. Oscillating test fields converge strongly two-scale:

\begin{proposition}\label{prop:admissibleproduct}
Let $W\in L^2(\Omega;C_\#(Y;\R^d))$. Then $W^\eps\twostrong W$, and so does every family $W^\eps+r_\eps$ with $r_\eps\to0$ in $L^2(\Omega;\R^d)$.
\end{proposition}

\begin{proof}
Choose $W_n\in C_c^\infty(\Omega;C_\#^\infty(Y;\R^d))$ with $W_n\to W$ in $L^2(\Omega;C_\#(Y;\R^d))$. For fixed $n$, periodic averaging gives $W_n^\eps\twoweak W_n$ and $\|W_n^\eps\|_{L^2(\Omega)}\to\|W_n\|_{L^2(\Omega\times Y)}$, while $\|(W-W_n)^\eps\|_{L^2(\Omega)}\le\|W-W_n\|_{L^2(\Omega;C_\#(Y;\R^d))}$ uniformly in $\eps$. Letting first $\eps\to0$ and then $n\to\infty$ gives $W^\eps\twostrong W$; the perturbation $r_\eps$ changes neither the two-scale limit nor the limit of the norms.
\end{proof}

The coefficient fields below are only measurable, so we also need the following multiplication property.

\begin{lemma}\label{lem:periodicmultiplier}
Let $B\in L_\#^\infty(Y;\R^{d\times d})$ and set $B_\eps(x)=B(x/\eps)$.
\begin{enumerate}[label=(\roman*)]
\item If $(z_\eps)_{\eps>0}$ is bounded in $L^2(\Omega;\R^d)$ and $z_\eps\twoweak z$, then $B_\eps z_\eps\twoweak B(y)z(x,y)$.
\item If $D\subset\R^N$ is measurable and $h\in L^2(D;\R^d)$, then $B_\eps h\weakto \overline B\,h$ in $L^2(D;\R^d)$, where $\overline B:=\int_YB(y)\,dy$.
\end{enumerate}
\end{lemma}

\begin{proof}
For smooth $B$ both assertions follow from the definition and periodic averaging. In general, let $B_n$ be periodic mollifications of $B$, so that $B_n\to B$ in $L^2(Y)$ and $\sup_n\|B_n\|_{L^\infty(Y)}\le\|B\|_{L^\infty(Y)}$. For a smooth test function $\varphi$, the error produced by $B-B_n$ is bounded by $\|z_\eps\|_{L^2(\Omega)}\|(B-B_n)^T(x/\eps)\varphi^\eps\|_{L^2(\Omega)}$, and periodic averaging gives
\[
 \lim_{\eps\to0}\bigl\|(B-B_n)^T(x/\eps)\varphi^\eps\bigr\|_{L^2(\Omega)}^2
 =\int_{\Omega\times Y}|(B-B_n)^T(y)\varphi(x,y)|^2\,dxdy\le C_\varphi\|B-B_n\|_{L^2(Y)}^2 .
\]
The limit of the error on the right is therefore bounded by
$C_\varphi^{1/2}\sup_\eps\|z_\eps\|_{L^2(\Omega)}\|B-B_n\|_{L^2(Y)}$, which tends
to zero as $n\to\infty$. The corresponding limiting integral is controlled
by the same $L^2(Y)$ approximation. This proves (i) first for smooth
periodic test functions $\varphi$ compactly supported in $x$.
For a general test $\varphi\in L^2(\Omega;C_\#(Y;\R^d))$,
approximate it in this Bochner norm by such smooth tests: the bound
$\|B_\eps z_\eps\|_{L^2(\Omega)}\le
\|B\|_\infty\sup_\eps\|z_\eps\|_{L^2(\Omega)}$,
together with $\|\varphi^\eps\|_{L^2(\Omega)}\le
\|\varphi\|_{L^2(\Omega;C_\#(Y))}$, makes the approximation
uniform in $\eps$, while $Bz\in L^2(\Omega\times Y)$ controls the
limit. For (ii), let $k\in L^2(D;\R^d)$; then $\int_DB_\eps h\cdot k\,dx=\int_{\R^N}B(x/\eps):(\mathbf 1_D\,k\otimes h)\,dx\to\int_D\overline B\,h\cdot k\,dx$, since $B(x/\eps)\overset{*}{\rightharpoonup}\overline B$ in $L^\infty(\R^N)$ by periodic averaging and $\mathbf 1_D\,k\otimes h\in L^1(\R^N)$.
\end{proof}

\section{Periodic fractional calculus}\label{sec:helmholtz}
\subsection{Periodic fractional Helmholtz decomposition}

The microscopic component of two-scale limits is a zero-mean periodic field. We first decompose such fields into fractional gradients and solenoidal fields, and then relate fractional and classical potentials. Define
\[
 \mathcal G_s:=\{\gradys v:v\in\Hperzero\},
 \qquad
 \mathcal Z_s:=\{Z\in L^2_{\#,0}(Y;\R^N):\divy Z=0\}.
\]

\begin{theorem}\label{thm:helmholtz}
For every $0<s<1$, $L^2_{\#,0}(Y;\R^N)=\mathcal G_s\overset{\perp}{\oplus}\mathcal Z_s$: every $F\in L^2_{\#,0}(Y;\R^N)$ decomposes uniquely as $F=\gradys v+Z$ with $v\in\Hperzero$ and $Z\in\mathcal Z_s$, and $\|F\|_{L^2(Y)}^2=\|\gradys v\|_{L^2(Y)}^2+\|Z\|_{L^2(Y)}^2$. The orthogonal projections are the classical Fourier Helmholtz projections; in particular, they do not depend on $s$.
\end{theorem}

\begin{proof}
For $k\ne0$, decompose
\[
 \widehat F(k)=P_k\widehat F(k)+(I-P_k)\widehat F(k),
 \qquad
 P_k:=\frac{q_k\otimes q_k}{|q_k|^2}.
\]
Set
\[
 \widehat Z(k)=(I-P_k)\widehat F(k),\qquad \widehat Z(0)=0.
\]
Then $q_k\cdot\widehat Z(k)=0$, hence $\divy Z=0$.
Indeed, for every zero-mean periodic vector field $W\in L^2_{\#,0}(Y;\R^N)$, the Fourier symbols yield
\[
 \divy W=0 \quad\Longleftrightarrow\quad
 q_k\cdot\widehat W(k)=0\quad(k\ne0)
 \quad\Longleftrightarrow\quad
 \operatorname{div}_y W=0
\]
in periodic distributions. Thus the solenoidal space $\mathcal Z_s$ is independent of $s$.

Define $v$ by
\begin{equation}\label{eq:potentialfourier}
 \widehat v(k)
 =\frac{q_k\cdot\widehat F(k)}{i|q_k|^{s+1}},\qquad k\ne0,
 \qquad \widehat v(0)=0.
\end{equation}
Then
\[
 iq_k|q_k|^{s-1}\widehat v(k)=P_k\widehat F(k),
\]
so $\gradys v=F-Z$. Moreover,
\[
 \|v\|_{\Hperzero}^2
 =\sum_{k\in\Z^N\setminus\{0\}}\frac{|q_k\cdot\widehat F(k)|^2}{|q_k|^2}
 \le \|F\|_{L^2(Y)}^2,
\]
which gives $v\in\Hperzero$. The range $\mathcal G_s$ is closed in $L^2_{\#,0}(Y;\R^N)$, since $\gradys$ is an isometry from the complete Hilbert space $\Hperzero$ onto $\mathcal G_s$. Orthogonality follows mode by mode from $P_k(I-P_k)=0$, and uniqueness follows from the orthogonal direct sum.
\end{proof}

The following de Rham-type consequence is the form used to identify microscopic limits.

\begin{corollary}\label{cor:deRham}
Let $F\in L^2_{\#,0}(Y;\R^N)$ satisfy $\int_YF\cdot\Phi\,dy=0$ for every $\Phi\in C_\#^\infty(Y;\R^N)$ with $\int_Y\Phi\,dy=0$ and $\divy\Phi=0$. Then $F=\gradys v$ for a unique $v\in\Hperzero$.
\end{corollary}

\begin{proof}
Smooth divergence-free trigonometric polynomials are dense in $\mathcal Z_s$ by Fourier truncation. Hence $F\perp\mathcal Z_s$, and Theorem~\ref{thm:helmholtz} applies.
\end{proof}

\subsection{The fractional--classical bridge and the two-scale state space}\label{sec:bridge}

Since the Helmholtz projections do not depend on $s$, the fractional and classical gradient spaces coincide. The following theorem makes this explicit at the level of potentials.

Define, on zero-mean periodic functions,
\[
 \mathcal J_s:=(-\Delta_y)^{(1-s)/2},
 \qquad
 \widehat{\mathcal J_s v}(k)=|q_k|^{1-s}\widehat v(k),\quad k\ne0.
\]
All powers of the periodic Laplacian in this subsection act on zero-mean
periodic functions; their zero Fourier mode is fixed to zero.

\begin{theorem}\label{thm:bridge}
The operator
\[
 \mathcal J_s:\dot H^1_\#(Y)\longrightarrow\dot H^s_\#(Y)
\]
is an onto isometry when both spaces are equipped with their homogeneous norms. Moreover,
\[
 \gradys(\mathcal J_s v)=\nabla_yv
 \qquad\forall v\in\dot H^1_\#(Y).
\]
Its inverse is $\mathcal J_s^{-1}=(-\Delta_y)^{-(1-s)/2}$ on zero-mean periodic functions. Consequently,
\begin{equation}\label{eq:rangeequality}
 \gradys\dot H^s_\#(Y)=\nabla_y\dot H^1_\#(Y)
 \quad\text{isometrically as closed subspaces of }L^2_{\#,0}(Y;\R^N).
\end{equation}
\end{theorem}

\begin{proof}
For $v\in\dot H^1_\#(Y)$,
\[
 \|\mathcal J_sv\|_{\dot H^s_\#}^2
 =\sum_{k\in\Z^N\setminus\{0\}}|q_k|^{2s}|q_k|^{2(1-s)}|\widehat v(k)|^2
 =\sum_{k\in\Z^N\setminus\{0\}}|q_k|^2|\widehat v(k)|^2
 =\|v\|_{\dot H^1_\#}^2.
\]
Furthermore,
\[
 \widehat{\gradys(\mathcal J_sv)}(k)
 =iq_k|q_k|^{s-1}|q_k|^{1-s}\widehat v(k)
 =iq_k\widehat v(k)
 =\widehat{\nabla_yv}(k).
\]
Given $w\in\dot H^s_\#(Y)$, define $v$ by
\[
 \widehat v(k)=|q_k|^{s-1}\widehat w(k),\qquad k\ne0.
\]
Then $v\in\dot H^1_\#(Y)$, $\mathcal J_sv=w$, and the inverse formula follows.
\end{proof}

As in the classical decomposition $\nabla u+\nabla_yu_1$ of two-scale limits of gradients \cite{Allaire1992}, \cite[Theorem~2.2]{VisintinFirstOrder2007}, we keep the macroscopic and microscopic potentials as separate components.

\begin{definition}\label{def:Htwo}
The two-scale state space is
\[
 \mathbb H_{2,0}^s(\Omega\times Y)
 :=H_0^s(\Omega)\times L^2\bigl(\Omega;\dot H_\#^s(Y)\bigr).
\]
For $U=(u_0,u_1)\in\mathbb H_{2,0}^s(\Omega\times Y)$, define the two-scale fractional gradient by
\begin{equation}\label{eq:twograddef}
 \twograds U(x,y):=\grads u_0(x)+\gradys u_1(x,y),\qquad(x,y)\in\Omega\times Y.
\end{equation}
A family $(u_\eps)_{\eps>0}$ \emph{realizes} $U$ if it is bounded in $H_0^s(\Omega)$, $u_\eps\to u_0$ in $L^2(\R^N)$, and $\grads u_\eps|_\Omega\twostrong\twograds U$ in $L^2(\Omega\times Y;\R^N)$.
\end{definition}

The zero-mean normalization of the microscopic potential makes the two contributions orthogonal.

\begin{proposition}\label{prop:twogradorth}
For every $U=(u_0,u_1)\in\mathbb H_{2,0}^s(\Omega\times Y)$,
\[
 \|\twograds U\|_{L^2(\Omega\times Y)}^2
 =\|\grads u_0\|_{L^2(\Omega)}^2
 +\|\gradys u_1\|_{L^2(\Omega\times Y)}^2.
\]
In particular, $\twograds U$ and $\grads u_0|_\Omega$ uniquely determine the zero-mean microscopic potential $u_1$.
\end{proposition}

\begin{proof}
The Fourier definition of $\gradys$ gives
\[
 \int_Y\gradys u_1(x,y)\,dy=0
 \quad\text{for a.e. }x.
\]
Hence the cross term in the square of \eqref{eq:twograddef} vanishes. Uniqueness follows from the injectivity of $\gradys$ on $\dot H_\#^s(Y)$.
\end{proof}

\section{Compactness and realization}\label{sec:oscdiv}
\subsection{Fractional two-scale compactness}\label{sec:compactness}

Compactness requires, besides the structure of the microscopic gradient space, a mechanism that eliminates fractional-divergence-free oscillations. It is provided by the following estimate.

\begin{lemma}\label{lem:oscdiv}
Let $\Phi\in C_c^\infty(\Omega;C_\#^\infty(Y;\R^N))$ satisfy
\begin{equation}\label{eq:solconditions}
 \int_Y\Phi(x,y)\,dy=0,
 \qquad
 \divy\Phi(x,y)=0\qquad(x\in\Omega,\ y\in Y),
\end{equation}
and let $\Phi_\eps(x)=\Phi(x,x/\eps)$ for $x\in\Omega$, extended by zero to $\R^N$. Then
\begin{equation}\label{eq:oscdivestimate}
 \|\divs\Phi_\eps\|_{L^2(\R^N)}
 \le C_\Phi\eps^{1-s};
\end{equation}
in particular, $\divs\Phi_\eps\to0$ strongly in $L^2(\R^N)$.
\end{lemma}

\begin{proof}
Expand
\[
 \Phi(x,y)=\sum_{k\in\Z^N\setminus\{0\}}\Phi_k(x)e^{iq_k\cdot y}.
\]
The condition $\divy\Phi=0$ gives
\[
 q_k\cdot\Phi_k(x)=0\qquad(k\ne0).
\]
It suffices to estimate one mode $F_{\eps,k}(x)=a(x)e^{iq_k\cdot x/\eps}$, $x\in\R^N$, with $q_k\cdot a(x)=0$. Plancherel gives
\[
 \|\divs F_{\eps,k}\|_2^2
 =\int_{\R^N}|(\xi+q_k/\eps)\cdot\widehat a(\xi)|^2
 |\xi+q_k/\eps|^{2(s-1)}\,d\xi.
\]
Since $q_k\cdot\widehat a(\xi)=0$, the scalar product equals $\xi\cdot\widehat a(\xi)$. On $E_1=\{|\xi|\le |q_k|/(2\eps)\}$ we have $|\xi+q_k/\eps|\ge|q_k|/(2\eps)$, so the integral over $E_1$ is at most $C|q_k/\eps|^{2(s-1)}\int_{\R^N}|\xi|^2|\widehat a(\xi)|^2\,d\xi$. On $E_2=\{|\xi|>|q_k|/(2\eps)\}$, using again $q_k\cdot\widehat a=0$, we have $|(\xi+q_k/\eps)\cdot\widehat a|\,|\xi+q_k/\eps|^{s-1}\le |\xi+q_k/\eps|^s|\widehat a|$ and $|\xi+q_k/\eps|^s\le C|\xi|^s\le C|q_k/\eps|^{s-1}|\xi|$, so the same estimate holds. Consequently,
\[
 \|\divs F_{\eps,k}\|_2
 \le C\eps^{1-s}|q_k|^{s-1}\|\nabla a\|_2.
\]
Since $\Phi$ is smooth and compactly supported in $x$, its Fourier coefficients decay rapidly in every $H^1(\R^N)$ seminorm. In particular, $\sum_{k\ne0}|q_k|^{s-1}\|\nabla_x\Phi_k\|_{L^2(\R^N)}<\infty$. Summing the modal estimates proves \eqref{eq:oscdivestimate} with $C_\Phi=C_{N,s}\sum_{k\ne0}|q_k|^{s-1}\|\nabla_x\Phi_k\|_2$, which depends on $N,s$ and the indicated seminorms of $\Phi$, but not on $\eps$.
\end{proof}

We now combine two-scale compactness in $L^2(\Omega)$, the oscillatory cancellation from Lemma~\ref{lem:oscdiv}, and the periodic Helmholtz characterization to identify the full two-scale limit of bounded fractional gradients.

\begin{theorem}\label{thm:maincompactness}
Let $(u_\eps)_{\eps>0}$ be uniformly bounded in $\Ho$. Given any null sequence $(\eps_n)_{n\in\mathbb N}$ with $\eps_n\downarrow0$, there exist a subsequence, not relabeled, and
\[
 u_0\in\Ho,
 \qquad
 u_1\in L^2(\Omega;\Hperzero)
\]
such that, as $n\to\infty$,
\begin{align}
 u_{\eps_n}&\weakto u_0 &&\text{in }H^s(\R^N),\label{eq:weakHs}\\
 u_{\eps_n}&\to u_0 &&\text{in }L^2(\R^N),\label{eq:strongL2}\\
 \grads u_{\eps_n}|_\Omega
 &\twoweak \grads u_0(x)+\gradys u_1(x,y)
 &&\text{in }L^2(\Omega\times Y;\R^N).\label{eq:main2scale}
\end{align}
Moreover,
\[
 \liminf_{n\to\infty}\|\grads u_{\eps_n}\|_{L^2(\Omega)}^2
 \ge
 \|\grads u_0\|_{L^2(\Omega)}^2
 +\|\gradys u_1\|_{L^2(\Omega\times Y)}^2.
\]
\end{theorem}

\begin{proof}
Fix the null sequence $(\eps_n)$ and, to lighten notation, write $\eps$ for $\eps_n$. By weak compactness in $H^s(\R^N)$ and compactness of $\Ho\hookrightarrow L^2(\Omega)$, after extraction we obtain \eqref{eq:weakHs}--\eqref{eq:strongL2}.

To identify the two-scale limit, observe that $\grads u_\eps$ is bounded in $L^2(\R^N;\R^N)$. Classical two-scale compactness on $\Omega$ gives
\[
 \grads u_\eps\twoweak G
 \quad\text{for some }G\in L^2(\Omega\times Y;\R^N).
\]
Since $\grads:H^s(\R^N)\to L^2(\R^N;\R^N)$ is linear and continuous,
\[
 \grads u_\eps\weakto \grads u_0
 \quad\text{in }L^2(\R^N;\R^N).
\]
Hence
\[
 \int_YG(x,y)\,dy=\grads u_0(x)
 \quad\text{for a.e. }x\in\Omega,
\]
and $G_1(x,y):=G(x,y)-\grads u_0(x)$ has zero $Y$-mean.

It remains to identify this zero-mean fluctuation. Let $\Phi$ satisfy \eqref{eq:solconditions}. By fractional integration by parts and Lemma~\ref{lem:oscdiv},
\[
 \left|\int_\Omega\grads u_\eps\cdot\Phi\!\left(x,\frac{x}{\eps}\right)dx\right|
 =\left|\int_{\R^N}u_\eps\divs\Phi_\eps\,dx\right|
 \le C\eps^{1-s}\to0.
\]
Passing to the two-scale limit gives
\[
 \int_{\Omega\times Y}G\cdot\Phi\,dxdy=0.
\]
The $y$-zero-mean condition on $\Phi$ removes the macroscopic term, so
\[
 \int_{\Omega\times Y}G_1\cdot\Phi\,dxdy=0.
\]
Testing with $\Phi(x,y)=\phi(x)\psi(y)$, $\phi\in C_c^\infty(\Omega)$, and using a countable dense family of smooth periodic solenoidal fields, one obtains for a.e.\ $x\in\Omega$
\[
 G_1(x,\cdot)\perp\mathcal Z_s.
\]
The orthogonality just obtained places $G_1$ in the potential component of the Helmholtz decomposition. Corollary~\ref{cor:deRham} therefore yields a unique $u_1(x,\cdot)\in\Hperzero$ such that
\[
 G_1(x,y)=\gradys u_1(x,y).
\]
The Fourier formula \eqref{eq:potentialfourier}, applied pointwise in $x$, gives measurability and
\[
 \int_\Omega\|u_1(x,\cdot)\|_{\Hperzero}^2dx
 =\|G_1\|_{L^2(\Omega\times Y)}^2.
\]
This proves \eqref{eq:main2scale}. Weak lower semicontinuity of the norm of $L^2(\Omega\times Y)$ under two-scale convergence then gives
\[
 \liminf_{\eps\to0}\|\grads u_\eps\|_{L^2(\Omega)}^2
 \ge \int_{\Omega\times Y}|\grads u_0+\gradys u_1|^2dxdy,
\]
and the cross term vanishes because $\int_Y\gradys u_1\,dy=0$.
\end{proof}

\begin{remark}\label{rem:oscdivs1}
For $s=1$ no rate is available: $\mathrm{div}\,\Phi_\eps=(\mathrm{div}_x\Phi)(x,x/\eps)$ is of order one, and the classical argument uses instead that $\int u_\eps(\mathrm{div}_x\Phi)^\eps\,dx\to\int_{\Omega\times Y}u_0\,\mathrm{div}_x\Phi\,dxdy=0$. For $s<1$ the smoothing $I_{1-s}$ of the oscillating zero-mean field $(\mathrm{div}_x\Phi)(x,x/\eps)$ produces the gain $\eps^{1-s}$, which degenerates as $s\uparrow1$.
\end{remark}

\begin{remark}\label{rem:transferproof}
An alternative route to the decomposition in Theorem~\ref{thm:maincompactness} uses classical two-scale compactness for gradients. Indeed, the Riesz potential $I_{1-s}$ maps $H_0^s(\Omega)$ continuously into $H^1_{\rm loc}(\R^N)$ and
\[
 \nabla I_{1-s}u=\grads u;
\]
see \cite{KreisbeckSchoenberger2022}. Applying classical gradient compactness to $I_{1-s}u_\eps$ on bounded sets containing $\Omega$ gives a microscopic classical gradient. By \eqref{eq:rangeequality}, that field is uniquely of the form $\gradys u_1$. The direct proof above gives in addition the rate $\eps^{1-s}$ in Lemma~\ref{lem:oscdiv}.
\end{remark}

\subsection{Realization at amplitude \texorpdfstring{$\eps^s$}{epsilon\string^s}}\label{sec:realization}

Theorem~\ref{thm:maincompactness} determines the form of every possible limit. Conversely, prescribed macroscopic and microscopic fields should be generated by admissible oscillating families. The realization rests on a quantitative estimate for oscillating functions with slowly varying amplitude, which in turn follows from a H\"older property of the symbol. Let $M_s(\xi)=i\xi|\xi|^{s-1}$, $\xi\in\R^N$, be the symbol of $\grads$; it is homogeneous of degree $s$.

\begin{lemma}\label{lem:symbolholder}
There exists $C_s>0$ such that
\[
 |M_s(\xi)-M_s(\zeta)|\le C_s|\xi-\zeta|^s
 \qquad(\xi,\zeta\in\R^N).
\]
\end{lemma}

\begin{proof}
Separate the cases $|\xi-\zeta|\ge\frac12\max\{|\xi|,|\zeta|\}$, where both terms are bounded by $C|\xi-\zeta|^s$, and its complement, where $|DM_s(z)|\le C_s|z|^{s-1}$ on the segment $[\xi,\zeta]$ gives the bound.
\end{proof}

For a single mode there is an exact identity: for $q\in2\pi\Z^N\setminus\{0\}$, $\grads(\eps^se^{iq\cdot x/\eps})=M_s(q)e^{iq\cdot x/\eps}$ in $\mathscr S'(\R^N)$. A microscopic potential of amplitude $\eps^s$ therefore generates an order-one fractional gradient, and Lemma~\ref{lem:amplitude} below is the stable version of this identity for a slowly varying amplitude.

\begin{lemma}\label{lem:amplitude}
Let $\omega(y)=\sum_{k\in\Z^N}c_ke^{iq_k\cdot y}$ with $\sum_{k\in\Z^N}(1+|q_k|^s)|c_k|<\infty$. Then, for every $a\in H^s(\R^N)$ and $\eps>0$,
\begin{equation}\label{eq:amplitudeestimate}
 \bigl\|\grads\bigl[\eps^sa(x)\omega(x/\eps)\bigr]-a(x)(\gradys\omega)(x/\eps)\bigr\|_{L^2(\R^N)}
 \le C_s\eps^s\|a\|_{\dot H^s(\R^N)}\sum_{k\in\Z^N}|c_k|,
\end{equation}
with $C_s$ the constant of Lemma~\ref{lem:symbolholder}.
\end{lemma}

\begin{proof}
Since $M_s(0)=0$ and $|M_s(q)|\le C_s|q|^s$, the series $\gradys\omega=\sum_{k\in\Z^N}M_s(q_k)c_ke^{iq_k\cdot y}$ converges uniformly. For one mode, the Fourier transform of $\grads[\eps^sa\,e^{iq_k\cdot x/\eps}]$ at $\xi+q_k/\eps$ is $\eps^sM_s(\xi+q_k/\eps)\widehat a(\xi)=M_s(q_k+\eps\xi)\widehat a(\xi)$, by homogeneity. Hence, by Plancherel and Lemma~\ref{lem:symbolholder},
\begin{equation}\label{eq:onemode}
\begin{aligned}
 \bigl\|\grads[\eps^sa\,e^{iq_k\cdot x/\eps}]-M_s(q_k)a\,e^{iq_k\cdot x/\eps}\bigr\|_{L^2(\R^N)}
 &=\bigl\|(M_s(q_k+\eps\xi)-M_s(q_k))\widehat a(\xi)\bigr\|_{L^2(\R^N_\xi)}\\
 &\le C_s\eps^s\|a\|_{\dot H^s(\R^N)} .
\end{aligned}
\end{equation}
Multiplying by $c_k$ and summing over $k\in\Z^N$ gives \eqref{eq:amplitudeestimate}.
\end{proof}

Applying this estimate to the Fourier series of the microscopic profile yields the realization result.

\begin{theorem}\label{thm:realization}
Let $u_0\in H_0^s(\Omega)$ and let $u_1\in C_c^\infty(\Omega;C_\#^\infty(Y))$ have zero $Y$-mean, with Fourier expansion $u_1(x,y)=\sum_{k\in\Z^N\setminus\{0\}}a_k(x)e^{iq_k\cdot y}$. Set
\[
 u_\eps(x)=u_0(x)+\eps^su_1\!\left(x,\frac{x}{\eps}\right),\qquad x\in\R^N.
\]
Then $u_\eps\in H_0^s(\Omega)$, $\sum_{k\in\Z^N\setminus\{0\}}\|a_k\|_{\dot H^s(\R^N)}<\infty$, and, for every $\eps>0$,
\begin{equation}\label{eq:strongrealization}
 \Bigl\|
 \grads u_\eps
 -\grads u_0
 -(\gradys u_1)\Bigl(x,\frac{x}{\eps}\Bigr)
 \Bigr\|_{L^2(\R^N)}
 \le C_s\,\eps^s\sum_{k\in\Z^N\setminus\{0\}}\|a_k\|_{\dot H^s(\R^N)}.
\end{equation}
Consequently, $(u_\eps)$ realizes $(u_0,u_1)$.
\end{theorem}

\begin{proof}
The oscillatory term is smooth and compactly supported in $\Omega$, so $u_\eps\in H_0^s(\Omega)$. Since $u_1$ is smooth in $y$ and compactly supported in $x$, the coefficients $a_k$ decay faster than every power of $|k|$ in every Sobolev norm, so the series in \eqref{eq:strongrealization} converges. The one-mode estimate \eqref{eq:onemode} with $a=a_k$, followed by summation over $k\in\Z^N\setminus\{0\}$, gives \eqref{eq:strongrealization}. Moreover, $\|u_\eps-u_0\|_{L^2(\R^N)}=\eps^s\|u_1(x,x/\eps)\|_{L^2(\R^N)}\le\eps^s\|u_1\|_{L^\infty(\Omega\times Y)}|\Omega|^{1/2}$. Finally, $W=\grads u_0+\gradys u_1$ belongs to $L^2(\Omega;C_\#(Y;\R^N))$, and Proposition~\ref{prop:admissibleproduct}, together with \eqref{eq:strongrealization}, gives the strong two-scale convergence.
\end{proof}

Theorem~\ref{thm:realization} uses the amplitude $\eps^s$. For the ansatz $u_0+\eps^\theta u_1(x,x/\eps)$, every other exponent gives either no microscopic component or an unbounded family.

\begin{proposition}\label{prop:amplitudedichotomy}
Let $u_0\in H_0^s(\Omega)$, let $u_1\in C_c^\infty(\Omega;C_\#^\infty(Y))$ have zero $Y$-mean with $u_1\not\equiv0$, let $\theta>0$, and set $u_\eps^\theta=u_0+\eps^\theta u_1(x,x/\eps)$.
\begin{enumerate}[label=(\roman*)]
\item If $\theta>s$, then $\grads u_\eps^\theta\to\grads u_0$ strongly in $L^2(\R^N)$; the two-scale limit has no microscopic component.
\item If $\theta<s$, then $\|\grads u_\eps^\theta\|_{L^2(\Omega)}\to\infty$; in particular $(u_\eps^\theta)$ is not bounded in $H^s(\R^N)$.
\end{enumerate}
\end{proposition}

\begin{proof}
By linearity, $\grads[\eps^\theta u_1(x,x/\eps)]=\eps^{\theta-s}\bigl[(\gradys u_1)(x,x/\eps)+R_\eps\bigr]$ with $\|R_\eps\|_{L^2(\R^N)}\le C\eps^s$ by \eqref{eq:strongrealization}. Periodic averaging gives $\|(\gradys u_1)(x,x/\eps)\|_{L^2(\Omega)}\to\|\gradys u_1\|_{L^2(\Omega\times Y)}$, which is positive because $\gradys$ is injective on zero-mean potentials and $u_1\not\equiv0$. For $\theta>s$ the prefactor $\eps^{\theta-s}$ tends to zero, which gives (i). For $\theta<s$,
\[
 \|\grads u_\eps^\theta\|_{L^2(\Omega)}
 \ge\eps^{\theta-s}\Bigl(\|(\gradys u_1)(x,x/\eps)\|_{L^2(\Omega)}-C\eps^s\Bigr)-\|\grads u_0\|_{L^2(\Omega)}\to\infty .
\]
\end{proof}

\subsection{Complete characterization of two-scale fractional gradients}\label{sec:characterization}

Compactness and realization together characterize the set of interior two-scale limits of Riesz fractional gradients.

\begin{theorem}\label{thm:characterization}
A field $G\in L^2(\Omega\times Y;\R^N)$ is the weak two-scale limit of $\grads u_\eps$ for some family $(u_\eps)_{\eps>0}$ bounded in $H_0^s(\Omega)$ if and only if $G=\grads u_0+\gradys u_1$ for some $u_0\in H_0^s(\Omega)$ and $u_1\in L^2(\Omega;\dot H^s_\#(Y))$. More precisely, every $(u_0,u_1)\in\mathbb H^s_{2,0}(\Omega\times Y)$ is realized by a family $(u_\eps)_{\eps>0}\subset H_0^s(\Omega)$; in particular,
\begin{equation}\label{eq:generalrecovery}
 u_\eps\to u_0\ \ \text{in }L^2(\R^N),
 \qquad
 \grads u_\eps\twostrong \grads u_0+\gradys u_1\ \ \text{in }L^2(\Omega\times Y;\R^N),
\end{equation}
and $u_1$ is uniquely determined by $G$ under the zero-mean normalization.
\end{theorem}

\begin{proof}
Necessity is Theorem~\ref{thm:maincompactness}. For the converse, we first approximate the prescribed microscopic potential by smooth finite-rank profiles. By density of smooth finite-rank functions in the Bochner space $L^2(\Omega;\dot H^s_\#(Y))$, choose $u_1^m\in C_c^\infty(\Omega;C_\#^\infty(Y))$ with zero $Y$-mean such that
\begin{equation}\label{eq:u1density}
 u_1^m\to u_1
 \quad\text{in }L^2(\Omega;\dot H^s_\#(Y)),
\end{equation}
and set $u_\eps^m=u_0+\eps^su_1^m(x,x/\eps)$. For fixed $m$, Theorem~\ref{thm:realization} yields $u_\eps^m\to u_0$ in $L^2(\R^N)$ and $\grads u_\eps^m\twostrong G_m:=\grads u_0+\gradys u_1^m$.
By \eqref{eq:u1density} and the periodic isometry,
\begin{equation}\label{eq:GmtoG}
 G_m\to G:=\grads u_0+\gradys u_1
 \quad\text{in }L^2(\Omega\times Y).
\end{equation}

It remains to choose the diagonal so that the resulting sequence is not only strongly two-scale convergent but also uniformly bounded in $H_0^s(\Omega)$. Fix a sequence $(\varphi_\ell)_{\ell\ge1}$ dense in $L^2(\Omega;C_\#(Y;\R^N))$. For each $m\ge1$, the strong realization theorem and periodic averaging allow us to choose $\eps_m>0$ recursively, with
$\eps_m<\min\{\eps_{m-1}/2,1/m\}$ for $m\ge2$ and $\eps_1<1$, so that for every $0<\eps\le\eps_m$,
\begin{align}
 &\|u_\eps^m-u_0\|_{L^2(\R^N)}\le 2^{-m},\label{eq:diagL2control}\\
 &\left\|\grads u_\eps^m-\grads u_0-(\gradys u_1^m)\!\left(x,\frac{x}{\eps}\right)\right\|_{L^2(\R^N)}\le 2^{-m},\label{eq:diagremcontrol}\\
 &\left|\left\|(\gradys u_1^m)\!\left(x,\frac{x}{\eps}\right)\right\|_{L^2(\Omega)}-
 \|\gradys u_1^m\|_{L^2(\Omega\times Y)}\right|\le2^{-m},\label{eq:diagnormcontrol}\\
 &\left|\|\grads u_\eps^m\|_{L^2(\Omega)}-
 \|G_m\|_{L^2(\Omega\times Y)}\right|\le2^{-m}.\label{eq:diagtotalnormcontrol}
\end{align}
In choosing $\eps_m$, we also impose, for every $1\le\ell\le m$ and every $0<\eps\le\eps_m$,
\begin{equation}\label{eq:diagtestcontrol}
 \left|\int_\Omega \grads u_\eps^m\cdot
 \varphi_\ell\!\left(x,\frac{x}{\eps}\right)\,dx
 -\int_{\Omega\times Y}G_m(x,y)\cdot\varphi_\ell(x,y)\,dxdy\right|
 \le 2^{-m}.
\end{equation}
This is possible because, for fixed $m$, $\grads u_\eps^m$ converges weakly two-scale to $G_m$, and there are only finitely many conditions at each stage. Set $m(\eps)=m$ if $\eps_{m+1}<\eps\le\eps_m$ and $u_\eps=u_\eps^{m(\eps)}$ for $0<\eps\le\eps_1$; for $\eps>\eps_1$, set $u_\eps=u_0$. This defines the family for every $\eps>0$. As $\eps\downarrow0$, $m(\eps)\to\infty$, and \eqref{eq:diagL2control} gives the first convergence in \eqref{eq:generalrecovery}. Moreover, \eqref{eq:diagremcontrol}--\eqref{eq:diagnormcontrol}, together with the fact that $(\gradys u_1^m)(x,x/\eps)$ vanishes outside $\Omega$, imply
\[
 \|\grads u_\eps\|_{L^2(\R^N)}
 \le \|\grads u_0\|_2
 +\|\gradys u_1^{m(\eps)}\|_{L^2(\Omega\times Y)}+2^{1-m(\eps)}.
\]
For $0<\eps\le\eps_1$, the right-hand side is uniformly bounded because $u_1^m\to u_1$ in $L^2(\Omega;\dot H_\#^s(Y))$. For $\eps>\eps_1$, we set $u_\eps=u_0$, so the same uniform bound holds on the remaining range. Together with the $L^2(\R^N)$ bound, this gives
\[
 \sup_{\eps>0}\|u_\eps\|_{H^s(\R^N)}<\infty.
\]
First, fix $\ell\ge1$. For all sufficiently small $\eps$, we have $m(\eps)\ge\ell$; by \eqref{eq:diagtestcontrol} and \eqref{eq:GmtoG},
\[
 \int_\Omega\grads u_\eps\cdot\varphi_\ell\!\left(x,\frac{x}{\eps}\right)\,dx
 \longrightarrow\int_{\Omega\times Y}G\cdot\varphi_\ell\,dxdy.
\]
The uniform $L^2$ bound on $\grads u_\eps$, the bound $\|\varphi^\eps\|_{L^2(\Omega)}\le\|\varphi\|_{L^2(\Omega;C_\#(Y))}$, and the density of $(\varphi_\ell)$ extend this convergence to all tests in $L^2(\Omega;C_\#(Y;\R^N))$. Thus
\[
 \grads u_\eps\twoweak G.
\]
Second, \eqref{eq:diagtotalnormcontrol} and \eqref{eq:GmtoG} give
\[
 \|\grads u_\eps\|_{L^2(\Omega)}\longrightarrow
 \|G\|_{L^2(\Omega\times Y)}.
\]
Combining weak two-scale convergence with convergence of the norms gives the strong two-scale convergence in \eqref{eq:generalrecovery}. Thus the diagonal family is bounded in $H_0^s(\Omega)$ and realizes the prescribed pair. Since $\int_YG\,dy=\grads u_0|_\Omega$, the uniqueness of $u_1$ follows from Proposition~\ref{prop:twogradorth}.
\end{proof}

\section{Elliptic homogenization and strong correctors}\label{sec:homogenization}
\subsection{Periodic problem and cell structure}\label{sec:cellstructure}

We first identify the homogenized problem (Subsections~\ref{sec:cellstructure}--\ref{sec:exterior}), then prove convergence of the energies in $\Omega$ and in $\Omega^c$ (Subsection~\ref{sec:energy}) and construct strong admissible correctors (Subsection~\ref{sec:correctors}). From now on,
\begin{equation}\label{eq:Aassumption}
 A\in L_\#^\infty(Y;\R^{N\times N}),
\end{equation}
with, for some $0<\alpha\le\beta$,
\begin{equation}\label{eq:ellipticity}
 A(y)\xi\cdot\xi\ge\alpha|\xi|^2,
 \qquad |A(y)\xi|\le\beta|\xi|
 \qquad\text{for all }\xi\in\R^N\text{ and a.e.\ }y\in Y.
\end{equation}
Symmetry of $A$ is not needed for the homogenization, energy and corrector results. It is used only in Lemma~\ref{lem:AhomAbar}, where we assume
\begin{equation}\label{eq:Asym}
 A(y)=A(y)^T\qquad\text{for a.e.\ }y\in Y.
\end{equation}
Set $A_\eps(x)=A(x/\eps)$, $x\in\R^N$. For $f\in(\Ho)'$, let $u_\eps\in\Ho$ solve
\begin{equation}\label{eq:microvariational}
 \int_{\R^N}A_\eps(x)\grads u_\eps\cdot\grads v\,dx
 =\langle f,v\rangle
 \qquad\forall v\in\Ho,
\end{equation}
that is, $-\divs(A_\eps\grads u_\eps)=f$ in $\Omega$ and $u_\eps=0$ in $\Omega^c$ in the variational sense. Taking $v=u_\eps$ and using \eqref{eq:ellipticity} and the fractional Poincar\'e inequality gives
\begin{equation}\label{eq:aprioribound}
 \|u_\eps\|_{H^s(\R^N)}\le C\|f\|_{(\Ho)'},
\end{equation}
so that Theorem~\ref{thm:maincompactness} applies.

For $\xi\in\R^N$, let $w_\xi^s\in\Hperzero$ be the unique solution, given by the Lax--Milgram theorem, of the cell problem
\begin{equation}\label{eq:cellfractional}
 \int_YA(y)\bigl(\xi+\gradys w_\xi^s\bigr)\cdot\gradys\varphi\,dy=0
 \qquad\forall\varphi\in\Hperzero.
\end{equation}
It depends linearly on $\xi$. Set $\chi_j^s:=w_{e_j}^s$ for $j=1,\dots,N$, $Q\in L^2_\#(Y;\R^{N\times N})$ the matrix with columns $\gradys\chi_1^s,\ldots,\gradys\chi_N^s$, and
\begin{equation}\label{eq:Ahomfractional}
 \Ahom\xi:=\int_YA(y)\bigl(\xi+\gradys w_\xi^s(y)\bigr)\,dy,\qquad \xi\in\R^N.
\end{equation}
Testing \eqref{eq:cellfractional} with $w_\xi^s$ gives
\begin{equation}\label{eq:cellenergy}
 \int_YA(y)\bigl(\xi+\gradys w_\xi^s\bigr)\cdot\bigl(\xi+\gradys w_\xi^s\bigr)\,dy=\Ahom\xi\cdot\xi,
 \qquad\xi\in\R^N,
\end{equation}
and hence $\Ahom\xi\cdot\xi\ge\alpha\|\xi+\gradys w_\xi^s\|_{L^2(Y)}^2\ge\alpha|\xi|^2$, since $\gradys w_\xi^s$ has zero mean. We also set $\Abar:=\int_YA(y)\,dy$ and introduce the bilinear forms of the microscopic and homogenized problems,
\begin{align*}
 \mathfrak a_\eps(u,v)
 &:=\int_{\R^N}A(x/\eps)\grads u\cdot\grads v\,dx,\\
 \mathfrak a_0(u,v)
 &:=\int_\Omega\Ahom\grads u\cdot\grads v\,dx
 +\int_{\Omega^c}\Abar\grads u\cdot\grads v\,dx,
\end{align*}
defined for $u,v\in H_0^s(\Omega)$. Thus \eqref{eq:microvariational} reads $\mathfrak a_\eps(u_\eps,v)=\langle f,v\rangle$, and the homogenized problem will be $\mathfrak a_0(u_0,v)=\langle f,v\rangle$.

\begin{proposition}\label{prop:classicalcell}
Let $\chi_j\in\dot H^1_\#(Y)$ be the classical cell correctors, that is, $\int_YA(y)(e_j+\nabla_y\chi_j)\cdot\nabla_y\varphi\,dy=0$ for all $\varphi\in\dot H^1_\#(Y)$. Then
\begin{equation}\label{eq:cellbridge}
 \chi_j^s=\mathcal J_s\chi_j=(-\Delta_y)^{(1-s)/2}\chi_j,
 \qquad
 \gradys\chi_j^s=\nabla_y\chi_j,
\end{equation}
and $\Ahom$ is the classical periodic homogenized matrix associated with $A$.
\end{proposition}

\begin{proof}
By Theorem~\ref{thm:bridge}, $\varphi\mapsto\mathcal J_s\varphi$ maps $\dot H^1_\#(Y)$ onto $\Hperzero$ with $\gradys\mathcal J_s\varphi=\nabla_y\varphi$. Hence \eqref{eq:cellfractional} is the classical cell problem written in the potential $\mathcal J_s^{-1}w_\xi^s$, and uniqueness gives \eqref{eq:cellbridge}. Inserting $\gradys\chi_j^s=\nabla_y\chi_j$ in \eqref{eq:Ahomfractional} gives the classical formula. In particular, the effective matrix $\Ahom$ is independent of $s\in(0,1)$, and the dependence on $s$ is carried entirely by the potentials $\chi_j^s=\mathcal J_s\chi_j$.
\end{proof}

For symmetric $A$, the cell potentials are nontrivial exactly when the two effective coefficients differ.

\begin{lemma}\label{lem:AhomAbar}
Let $A$ satisfy \eqref{eq:Aassumption}--\eqref{eq:Asym}. Then $\Ahom=\Abar$ if and only if $Q\equiv0$, that is, if and only if $\mathrm{div}_y(A(y)e_j)=0$ in $\mathcal D'(\T^N)$ for $j=1,\dots,N$.
\end{lemma}

\begin{proof}
By Proposition~\ref{prop:classicalcell} and the minimum principle for the classical cell problem, $\Ahom\xi\cdot\xi=\min_{\varphi\in\dot H^1_\#(Y)}\int_YA(\xi+\nabla_y\varphi)\cdot(\xi+\nabla_y\varphi)\,dy\le\Abar\xi\cdot\xi$. By strict convexity, equality holds for all $\xi\in\R^N$ if and only if $\varphi=0$ is the minimizer for every $\xi\in\R^N$, i.e.\ if and only if $\int_YA\xi\cdot\nabla_y\varphi\,dy=0$ for all $\varphi\in\dot H^1_\#(Y)$ and $\xi\in\R^N$. For symmetric $A$ both $\Ahom$ and $\Abar$ are symmetric, so equality of the quadratic forms is equality of the matrices.
\end{proof}

\subsection{Exterior pairing and homogenization}\label{sec:exterior}

Unlike the classical case, the bilinear form $\mathfrak a_\eps$ also involves $\Omega^c$, where $\grads u_\eps$ does not vanish. Its exterior part is controlled by a far-field estimate for the kernel of $\grads$. Recall that $I_{1-s}$ has the convolution kernel $c_{N,1-s}|x|^{1-s-N}$, where
\[
 c_{N,\varsigma}=\frac{\Gamma\bigl(\frac{N-\varsigma}{2}\bigr)}{2^{\varsigma}\pi^{N/2}\Gamma\bigl(\frac{\varsigma}{2}\bigr)},\qquad 0<\varsigma<N .
\]
Differentiating, for $z\in C_c^\infty(\R^N)$,
\begin{equation}\label{eq:gradskernel}
\begin{gathered}
 \grads z(x)=\mathrm{p.v.}\!\int_{\R^N}K_s(x-y)\,z(y)\,dy,\\
 K_s(w)=(1-s-N)\,c_{N,1-s}\,\frac{w}{|w|^{N+1+s}}\quad(w\in\R^N\setminus\{0\}),\\
 |K_s(w)|\le C_{N,s}|w|^{-N-s}.
\end{gathered}
\end{equation}
The principal-value representation in \eqref{eq:gradskernel} is initially
understood for $z\in C_c^\infty(\R^N)$; the corresponding convolution
identity extends in the sense of distributions whenever the terms are
defined. In particular, away from $\operatorname{supp}z$, it reduces to
an ordinary absolutely convergent integral for compactly supported
$z\in L^2(\R^N)$.

\begin{lemma}\label{lem:farfield}
Let $\mathcal K\subset\R^N$ be compact, let $\Sigma\subset\R^N$ be measurable with $\dist(\Sigma,\mathcal K)\ge d>0$, and let $z\in L^2(\R^N)$ vanish outside $\mathcal K$. Then $\grads z$ coincides on $\Sigma$ with an absolutely convergent integral and
\begin{equation}\label{eq:farfieldgeneral}
 \|\grads z\|_{L^2(\Sigma)}\le C_{N,s}\,d^{-s}\,\|z\|_{L^2(\R^N)} .
\end{equation}
\end{lemma}

\begin{proof}
The kernel $K_s$ is locally integrable away from the origin.
For compactly supported $z\in L^2(\R^N)$, its distributional fractional
gradient agrees, on every open set disjoint from $\mathcal K$, with the ordinary
convolution $K_s*z$; this follows by testing the distributional convolution
against compactly supported smooth functions on that open set and using
Fubini, since the kernel has no singularity on the separated supports.
Thus for almost every $x\in\Sigma$ the integral over $\mathcal K$ is absolutely
convergent (indeed $z\in L^1(\mathcal K)$), and the same integral equals
$(K_s\mathbf 1_{\{|w|\ge d\}})*z(x)$. Young's inequality and
$\|K_s\mathbf 1_{\{|w|\ge d\}}\|_{L^1(\R^N)}\le
 C_{N,s}d^{-s}$ give \eqref{eq:farfieldgeneral}.
\end{proof}

If $\Omega\Subset\Omega_1$, then $u\mapsto\grads u|_{\Omega_1^c}$ is compact from $H_0^s(\Omega)$ into $L^2(\Omega_1^c;\R^N)$; see \cite[Lemma~2.12]{KreisbeckSchoenberger2022} for the complementary-value setting and \cite[Proposition~2.10]{BruscaEtAl2026} for the formulation as strong convergence locally in $\R^N\setminus\overline\Omega$. This also follows from Lemma~\ref{lem:farfield} and the compact embedding $H_0^s(\Omega)\Subset L^2(\Omega)$. Such estimates control the gradient away from $\partial\Omega$; for the homogenized equation, the following exterior pairing suffices.

\begin{lemma}\label{lem:exteriorpairing}
Assume $u_\eps\weakto u$ in $H_0^s(\Omega)$ and let
\[
 B\in L_\#^\infty(Y;\R^{N\times N}),
 \qquad B_\eps(x)=B(x/\eps),
 \qquad \overline B=\int_YB(y)\,dy.
\]
Then, for every $h\in L^2(\Omega^c;\R^N)$,
\begin{equation}\label{eq:exteriorpairing}
 \int_{\Omega^c}B_\eps\grads u_\eps\cdot h\,dx
 \longrightarrow
 \int_{\Omega^c}\overline B\,\grads u\cdot h\,dx.
\end{equation}
\end{lemma}

\begin{proof}
The family $(\grads u_\eps)_{\eps>0}$ is uniformly bounded in $L^2(\R^N)$. Let $\varpi>0$. The sets $\Omega_t\setminus\Omega$, $\Omega_t:=\Omega+\mathbb B_t$, decrease to $\overline\Omega\setminus\Omega=\partial\Omega$ as $t\downarrow0$. Since $|\partial\Omega|=0$, dominated convergence gives a bounded open set $\Omega_1=\Omega_t$ with $\Omega\Subset\Omega_1$ such that
\[
 \|h\|_{L^2(\Omega_1\setminus\Omega)}<\varpi.
\]
The compact embedding $H_0^s(\Omega)\Subset L^2(\Omega)$ gives
$u_\eps\to u$ strongly in $L^2(\Omega)$. Since the difference
$u_\eps-u$ is supported in $\overline\Omega$, Lemma~\ref{lem:farfield}
therefore yields the quantitative estimate
\[
 \|\grads(u_\eps-u)\|_{L^2(\Omega_1^c)}
 \le C_{\Omega,\Omega_1,s}\|u_\eps-u\|_{L^2(\Omega)}\longrightarrow0.
\]
By Lemma~\ref{lem:periodicmultiplier}(ii) with $D=\Omega_1^c$, $B_\eps^Th\weakto\overline B^{\,T}h$ in $L^2(\Omega_1^c;\R^N)$, and therefore
\[
 \int_{\Omega_1^c}B_\eps\grads u_\eps\cdot h\,dx
 \to
 \int_{\Omega_1^c}\overline B\,\grads u\cdot h\,dx.
\]
For the remaining layer, uniform boundedness gives
\[
 \left|\int_{\Omega_1\setminus\Omega}
 B_\eps\grads u_\eps\cdot h\,dx\right|
 \le \|B\|_{L^\infty(Y)}
 \|\grads u_\eps\|_2
 \|h\|_{L^2(\Omega_1\setminus\Omega)}
 \le C\varpi,
\]
and the same estimate holds for the limiting integral. Letting $\varpi\downarrow0$ proves \eqref{eq:exteriorpairing}.
\end{proof}

In $\Omega$ the oscillations produce $\Ahom$; in $\Omega^c$, Lemma~\ref{lem:exteriorpairing} with $B=A$ shows that the surviving coefficient is $\Abar$.

\begin{theorem}\label{thm:homogenization}
Let \eqref{eq:Aassumption}--\eqref{eq:ellipticity} hold, and let $u_\eps$ solve \eqref{eq:microvariational}. Then the whole family converges, weakly in $\Ho$ and strongly in $L^2(\R^N)$, to the unique $u_0\in\Ho$ satisfying
\begin{equation}\label{eq:homogenizedvariational}
 \int_\Omega\Ahom\grads u_0\cdot\grads v\,dx
 +\int_{\Omega^c}\Abar\grads u_0\cdot\grads v\,dx
 =\langle f,v\rangle
 \qquad\forall v\in\Ho.
\end{equation}
Moreover, the microscopic potential of Theorem~\ref{thm:maincompactness} is
\begin{equation}\label{eq:u1formula}
 u_1(x,y)=\sum_{j=1}^N(\grads u_0(x))_j\chi_j^s(y),\qquad(x,y)\in\Omega\times Y,
\end{equation}
and the whole family satisfies
\begin{equation}\label{eq:gradlimitcell}
 \grads u_\eps\twoweak
 \left[\Id+Q(y)\right]\grads u_0(x)
 \quad\text{in }L^2(\Omega\times Y;\R^N).
\end{equation}
\end{theorem}

\begin{proof}
The a priori bound \eqref{eq:aprioribound} and Theorem~\ref{thm:maincompactness} give, along a subsequence,
\[
 u_\eps\weakto u_0\quad\text{in }H_0^s(\Omega),
 \qquad
 \grads u_\eps\twoweak G:=\grads u_0+\gradys u_1
 \quad\text{in }L^2(\Omega\times Y;\R^N).
\]
We first identify the interior flux. By Lemma~\ref{lem:periodicmultiplier}(i),
\begin{equation}\label{eq:flux2scale}
 A_\eps\grads u_\eps\twoweak A(y)G(x,y)
 \quad\text{in }L^2(\Omega\times Y;\R^N).
\end{equation}

To identify both the macroscopic and microscopic equations, we use realized oscillatory test functions: let $v_0\in C_c^\infty(\Omega)$, let $v_1\in C_c^\infty(\Omega;C_\#^\infty(Y))$ have zero $Y$-mean, and set $v_\eps=v_0+\eps^sv_1(x,x/\eps)$. Then $v_\eps\in H_0^s(\Omega)$, and Theorem~\ref{thm:realization} gives
\begin{equation}\label{eq:testgradientrealization}
 \grads v_\eps-
 \left[\grads v_0+(\gradys v_1)\!\left(x,\frac{x}{\eps}\right)\right]
 \to0
 \quad\text{in }L^2(\R^N).
\end{equation}
In particular, $(v_\eps)_{\eps>0}$ is bounded in $H_0^s(\Omega)$ and converges to $v_0$ in $L^2(\Omega)$, so that $v_\eps\weakto v_0$ in $H_0^s(\Omega)$ and $\langle f,v_\eps\rangle\to\langle f,v_0\rangle$.

We split the variational identity into $\Omega$ and $\Omega^c$. On $\Omega$, $\grads v_\eps=\Psi^\eps+r_\eps$ with $\Psi:=\grads v_0+\gradys v_1\in L^2(\Omega;C_\#(Y;\R^N))$ and $r_\eps\to0$ in $L^2(\Omega)$ by \eqref{eq:testgradientrealization}; testing \eqref{eq:flux2scale} with $\Psi$ gives
\begin{equation}\label{eq:interior2scalelimit}
 \int_\Omega A_\eps\grads u_\eps\cdot\grads v_\eps\,dx
 \longrightarrow
 \int_{\Omega\times Y}
 A(y)G(x,y)\cdot
 \bigl(\grads v_0(x)+\gradys v_1(x,y)\bigr)\,dxdy.
\end{equation}
Outside $\Omega$, the oscillatory part of $v_\eps$ has fractional gradient converging strongly to zero by \eqref{eq:testgradientrealization}; hence
\[
 \int_{\Omega^c}A_\eps\grads u_\eps\cdot\grads v_\eps\,dx
 -
 \int_{\Omega^c}A_\eps\grads u_\eps\cdot\grads v_0\,dx
 \to0.
\]
Applying Lemma~\ref{lem:exteriorpairing} with $B=A$ and $h=\grads v_0|_{\Omega^c}$ yields
\[
 \int_{\Omega^c}A_\eps\grads u_\eps\cdot\grads v_\eps\,dx
 \to
 \int_{\Omega^c}\Abar\grads u_0\cdot\grads v_0\,dx.
\]
Combining the interior, exterior, and right-hand-side limits, we obtain, for all such $(v_0,v_1)$, the two-scale variational identity
\begin{equation}\label{eq:twoscalevariational}
 \int_{\Omega\times Y}
 A(y)(\grads u_0+\gradys u_1)
 \cdot(\grads v_0+\gradys v_1)\,dxdy
 +\int_{\Omega^c}\Abar\grads u_0\cdot\grads v_0\,dx
 =\langle f,v_0\rangle.
\end{equation}

To identify the microscopic field, set $v_0=0$ in \eqref{eq:twoscalevariational} and take $v_1(x,y)=\phi(x)\varphi(y)$ with $\phi\in C_c^\infty(\Omega)$ and $\varphi\in C_\#^\infty(Y)$ of zero mean. Localization in $x$ gives, for a.e.\ $x\in\Omega$,
\[
 \int_YA(y)
 \bigl(\grads u_0(x)+\gradys u_1(x,y)\bigr)
 \cdot\gradys\varphi(y)\,dy=0.
\]
By density in $\Hperzero$ and uniqueness of the cell problem \eqref{eq:cellfractional},
\[
 u_1(x,\cdot)=w_{\grads u_0(x)}^s
 =\sum_{j=1}^N(\grads u_0(x))_j\chi_j^s,
\]
which proves \eqref{eq:u1formula} and \eqref{eq:gradlimitcell}.

The macroscopic equation follows by setting $v_1=0$ in \eqref{eq:twoscalevariational}. By the definition of $\Ahom$,
\[
 \int_YA(y)
 \bigl(\grads u_0+\gradys u_1\bigr)\,dy
 =\Ahom\grads u_0,
\]
so
\[
 \int_\Omega\Ahom\grads u_0\cdot\grads v_0\,dx
 +\int_{\Omega^c}\Abar\grads u_0\cdot\grads v_0\,dx
 =\langle f,v_0\rangle
\]
for every $v_0\in C_c^\infty(\Omega)$. By the definition \eqref{eq:H0s} and continuity of the bilinear form, this extends to every $v\in H_0^s(\Omega)$, proving \eqref{eq:homogenizedvariational}.
By \eqref{eq:cellenergy}, $\Ahom\xi\cdot\xi\ge\alpha|\xi|^2$; likewise $\Abar\xi\cdot\xi\ge\alpha|\xi|^2$.
Consequently the homogenized form is bounded and coercive on $\Ho$,
by the fractional Poincar\'e inequality. Its solution is unique, so the
whole family, not merely the extracted subsequence, converges.
\end{proof}

\begin{remark}\label{rem:effectiveknown}
The coefficient law in \eqref{eq:homogenizedvariational} agrees with \cite[Theorem~1.3]{BuchingerEtAl2026}: local periodic $H$-convergence gives $\Ahom$ in $\Omega$, while weak-* averaging gives $\Abar$ outside $\overline\Omega$. Theorem~\ref{thm:homogenization} adds the microscopic identification \eqref{eq:u1formula}--\eqref{eq:gradlimitcell}, on which the correctors of Subsection~\ref{sec:correctors} are built.
\end{remark}

\begin{example}\label{ex:laminatejump}
For symmetric $A$, Lemma~\ref{lem:AhomAbar} shows that the coefficient of \eqref{eq:homogenizedvariational} jumps across $\partial\Omega$ exactly when $Q\ne0$, that is, when the cell potentials are nontrivial. For instance, let $A(y)=a(y_1)\Id$ be a laminate, with $a\in L^\infty_\#(0,1)$, $\alpha\le a\le\beta$, and write $\overline b:=\int_0^1b(t)\,dt$ for the mean of a $1$-periodic function $b$. The classical cell problems give $\chi_j=0$ for $j\ge2$ and $\partial_{y_1}\chi_1=(\overline{a^{-1}})^{-1}a^{-1}-1$. Hence
\[
 \Ahom=\operatorname{diag}\bigl((\overline{a^{-1}})^{-1},\overline a,\dots,\overline a\bigr),
 \qquad
 \Abar=\overline a\,\Id,
\]
and $Q$ has the single nonzero column $(\partial_{y_1}\chi_1)e_1$. By the Cauchy--Schwarz inequality, $(\overline{a^{-1}})^{-1}\le\overline a$, with equality only if $a$ is constant. As soon as $a$ is not constant, the coefficient of \eqref{eq:homogenizedvariational} therefore jumps across $\partial\Omega$ in the direction $e_1$: the harmonic mean acts inside $\Omega$ and the arithmetic mean outside.
\end{example}

\subsection{Energy convergence}\label{sec:energy}

The variational structure upgrades the weak two-scale convergence of Theorem~\ref{thm:homogenization} to convergence of the energies, separately in $\Omega$ and in $\Omega^c$. Only the symmetric part
\[
 A^{\rm sym}:=\tfrac12\bigl(A+A^T\bigr)
\]
of $A$ enters a quadratic form; it satisfies \eqref{eq:ellipticity} with the same $\alpha$ and $\beta$, since $A^{\rm sym}\xi\cdot\xi=A\xi\cdot\xi$ and $|A^{\rm sym}\xi|\le\frac12(|A\xi|+|A^T\xi|)\le\beta|\xi|$, and it coincides with $A$ when $A$ is symmetric.

\begin{theorem}\label{thm:strongenergy}
Under the assumptions of Theorem~\ref{thm:homogenization}, let
\[
 G(x,y):=\bigl[\Id+Q(y)\bigr]\grads u_0(x),\qquad(x,y)\in\Omega\times Y.
\]
Then
\begin{equation}\label{eq:interiorenergyconv}
 \int_\Omega A(x/\eps)\grads u_\eps\cdot\grads u_\eps\,dx
 \longrightarrow
 \int_{\Omega\times Y}A(y)G\cdot G\,dxdy
 =\int_\Omega\Ahom\grads u_0\cdot\grads u_0\,dx
\end{equation}
and
\begin{equation}\label{eq:exteriorenergyconv}
 \int_{\Omega^c}A(x/\eps)\grads u_\eps\cdot\grads u_\eps\,dx
 \longrightarrow
 \int_{\Omega^c}\Abar\grads u_0\cdot\grads u_0\,dx .
\end{equation}
Moreover, with $(A^{\rm sym})^{1/2}$ the positive symmetric square root of $A^{\rm sym}$,
\begin{equation}\label{eq:strongenergy}
 (A^{\rm sym})^{1/2}(x/\eps)\grads u_\eps
 \twostrong
 (A^{\rm sym})^{1/2}(y)G(x,y)
 \qquad\text{in }L^2(\Omega\times Y;\R^N).
\end{equation}
For symmetric $A$, \eqref{eq:strongenergy} reads $A^{1/2}(x/\eps)\grads u_\eps\twostrong A^{1/2}(y)G$.
\end{theorem}

\begin{proof}
Write $g_\eps:=\grads u_\eps$ and $S:=A^{\rm sym}$. With $\xi=\grads u_0(x)$, \eqref{eq:cellenergy} and \eqref{eq:gradlimitcell} yield the identity in \eqref{eq:interiorenergyconv}. Since $Av\cdot v=Sv\cdot v$ for every $v\in\R^N$, all quadratic expressions below may be written with $S$.

Testing \eqref{eq:microvariational} with $u_\eps$ and \eqref{eq:homogenizedvariational} with $u_0$ gives convergence of the total energies,
\[
 \mathfrak a_\eps(u_\eps,u_\eps)=\langle f,u_\eps\rangle\longrightarrow\langle f,u_0\rangle=\mathfrak a_0(u_0,u_0),
\]
because $u_\eps\weakto u_0$ in $\Ho$.

For the interior lower bound, let $\Psi\in L^2(\Omega;C_\#(Y;\R^N))$. Since $S$ is positive semidefinite,
\[
 0\le\int_\Omega S_\eps(g_\eps-\Psi^\eps)\cdot(g_\eps-\Psi^\eps)\,dx
 =\int_\Omega S_\eps g_\eps\cdot g_\eps\,dx-2\int_\Omega S_\eps g_\eps\cdot\Psi^\eps\,dx+\int_\Omega S_\eps\Psi^\eps\cdot\Psi^\eps\,dx .
\]
By Lemma~\ref{lem:periodicmultiplier}(i), applied to $g_\eps\twoweak G$ and to $\Psi^\eps\twoweak\Psi$ (Proposition~\ref{prop:admissibleproduct}), $S_\eps g_\eps\twoweak SG$ and $S_\eps\Psi^\eps\twoweak S\Psi$; testing with $\Psi$ gives
\[
 \liminf_{\eps\to0}\int_\Omega S_\eps g_\eps\cdot g_\eps\,dx
 \ge\int_{\Omega\times Y}\bigl(2SG\cdot\Psi-S\Psi\cdot\Psi\bigr)\,dxdy .
\]
Choosing $\Psi_n\to G$ in $L^2(\Omega\times Y;\R^N)$, which is possible by density, we obtain
\[
 \liminf_{\eps\to0}\int_\Omega A_\eps g_\eps\cdot g_\eps\,dx
 \ge\int_{\Omega\times Y}SG\cdot G\,dxdy=\int_{\Omega\times Y}AG\cdot G\,dxdy .
\]

For the exterior lower bound, let $t>0$ and $\Omega_1:=\Omega_t=\Omega+\mathbb B_t$. Since $\dist(\Omega_1^c,\overline\Omega)\ge t$ and $u_\eps\to u_0$ in $L^2(\R^N)$, Lemma~\ref{lem:farfield} with $\mathcal K=\overline\Omega$ and $d=t$ gives $g_\eps\to\grads u_0$ in $L^2(\Omega_1^c;\R^N)$, and Lemma~\ref{lem:periodicmultiplier}(ii) with $D=\Omega_1^c$ and $h=\grads u_0$ gives $\int_{\Omega_1^c}A_\eps\grads u_0\cdot\grads u_0\to\int_{\Omega_1^c}\Abar\grads u_0\cdot\grads u_0$; hence
\[
 \int_{\Omega_1^c}A_\eps g_\eps\cdot g_\eps\,dx\longrightarrow\int_{\Omega_1^c}\Abar\grads u_0\cdot\grads u_0\,dx .
\]
Since $A_\eps g_\eps\cdot g_\eps\ge0$ on $\Omega_1\setminus\Omega$, the right-hand side bounds $\liminf_{\eps\to0}\int_{\Omega^c}A_\eps g_\eps\cdot g_\eps\,dx$ from below; letting $t\downarrow0$, so that $\Omega_t^c\uparrow\R^N\setminus\overline\Omega$, and using $|\partial\Omega|=0$, we obtain
\[
 \liminf_{\eps\to0}\int_{\Omega^c}A_\eps g_\eps\cdot g_\eps\,dx
 \ge\int_{\Omega^c}\Abar\grads u_0\cdot\grads u_0\,dx .
\]

Write $E_\eps^{\rm int}$ and $E_\eps^{\rm ext}$ for the interior and exterior energies of $u_\eps$, and $E_0^{\rm int}$ and $E_0^{\rm ext}$ for the right-hand sides of the two lower bounds, so that $E_0^{\rm int}+E_0^{\rm ext}=\mathfrak a_0(u_0,u_0)$. Since $E_\eps^{\rm int}+E_\eps^{\rm ext}\to E_0^{\rm int}+E_0^{\rm ext}$,
\[
 \limsup_{\eps\to0}E_\eps^{\rm int}
 \le E_0^{\rm int}+E_0^{\rm ext}-\liminf_{\eps\to0}E_\eps^{\rm ext}
 \le E_0^{\rm int},
\]
and symmetrically $\limsup_{\eps\to0}E_\eps^{\rm ext}\le E_0^{\rm ext}$. Together with the lower bounds, this proves \eqref{eq:interiorenergyconv} and \eqref{eq:exteriorenergyconv}. Finally, $(S^{1/2})_\eps g_\eps\twoweak S^{1/2}G$ by Lemma~\ref{lem:periodicmultiplier}(i), and
\[
 \|(S^{1/2})_\eps g_\eps\|_{L^2(\Omega)}^2=\int_\Omega A_\eps g_\eps\cdot g_\eps\,dx
 \longrightarrow\int_{\Omega\times Y}SG\cdot G\,dxdy=\|S^{1/2}G\|_{L^2(\Omega\times Y)}^2,
\]
which gives \eqref{eq:strongenergy}.
\end{proof}

\subsection{Strong correctors}\label{sec:correctors}

We now construct strong correctors for the fractional gradients and then
for the solutions. Set $p:=\grads u_0\in L^2(\R^N;\R^N)$.
Since $A\in L^\infty_\#(Y;\R^{N\times N})$ guarantees only
$Q\in L^2_\#(Y;\R^{N\times N})$, the formal corrected gradient
$[\Id+Q(x/\eps)]p(x)$ must be approximated.

For $j=1,\dots,N$ choose zero-mean trigonometric polynomials $\chi_{j,m}^s$, $m\in\N$, with
\[
 \chi_{j,m}^s\to\chi_j^s
 \quad\text{in }\Hperzero\quad(m\to\infty),
\]
and let $Q_m\in L^\infty_\#(Y;\R^{N\times N})$ be the matrix with columns $\gradys\chi_{1,m}^s,\ldots,\gradys\chi_{N,m}^s$. By \eqref{eq:isometry}, $Q_m\to Q$ in $L^2(Y;\R^{N\times N})$. Let $\rho_r$ be a standard mollifier and $p_r:=\rho_r*p$, so that $p_r\to p$ in $L^2(\R^N;\R^N)$ as $r\downarrow0$. Define
\[
 G_\eps^{m,r}(x)
 :=p(x)+\mathbf1_\Omega(x)Q_m(x/\eps)p_r(x),\qquad x\in\R^N.
\]

\begin{proposition}\label{prop:gradientcorrector}
For fixed $m\in\N$ and $r>0$,
\begin{equation}\label{eq:energycorrector}
 \limsup_{\eps\to0}
 \|\grads u_\eps-G_\eps^{m,r}\|_{L^2(\R^N)}^2
 \le
 \frac{\beta}{\alpha}
 \|Q_mp_r-Qp\|_{L^2(\Omega\times Y)}^2.
\end{equation}
Consequently,
\begin{equation}\label{eq:doublelimitcorrector}
 \lim_{r\downarrow0}\lim_{m\to\infty}\limsup_{\eps\to0}
 \|\grads u_\eps-G_\eps^{m,r}\|_{L^2(\R^N)}=0.
\end{equation}
\end{proposition}

\begin{proof}
By \eqref{eq:ellipticity},
\begin{equation}\label{eq:Eepsmr}
 \alpha\|\grads u_\eps-G_\eps^{m,r}\|_{L^2(\R^N)}^2\le E_\eps:=\int_{\R^N}A_\eps(\grads u_\eps-G_\eps^{m,r})\cdot(\grads u_\eps-G_\eps^{m,r})\,dx .
\end{equation}
Put $G:=(\Id+Q)p$ and $H:=p+Q_mp_r$ on $\Omega\times Y$. We pass to the limit in the four terms of the expansion of $E_\eps$. By Theorem~\ref{thm:strongenergy}, the quadratic term tends to $\int_{\Omega\times Y}AG\cdot G\,dxdy+\int_{\Omega^c}\Abar p\cdot p\,dx$. In $\Omega$, $A_\eps\grads u_\eps\twoweak AG$ and $A_\eps^T\grads u_\eps\twoweak A^TG$ by Lemma~\ref{lem:periodicmultiplier}(i) and Theorem~\ref{thm:homogenization}, while on $\Omega$ we have $G_\eps^{m,r}=H^\eps$ with $H=p+Q_mp_r\in L^2(\Omega;C_\#(Y;\R^N))$, and $A_\eps H^\eps\twoweak AH$ by Lemma~\ref{lem:periodicmultiplier}(i) and Proposition~\ref{prop:admissibleproduct}; testing these three weak two-scale limits with $H$ gives the interior limits of the mixed and quadratic terms in $G_\eps^{m,r}$. In $\Omega^c$ we have $G_\eps^{m,r}=p$; Lemma~\ref{lem:exteriorpairing} with $B=A$, $B=A^T$ and $h=p|_{\Omega^c}$, and Lemma~\ref{lem:periodicmultiplier}(ii) with $D=\Omega^c$, show that the three remaining exterior contributions all tend to $\int_{\Omega^c}\Abar p\cdot p\,dx$, so that the exterior terms cancel in the limit. Since $A^TG\cdot H=AH\cdot G$, the interior limit is
\[
 \int_{\Omega\times Y}\bigl(AG\cdot G-AG\cdot H-A^TG\cdot H+AH\cdot H\bigr)\,dxdy
 =\int_{\Omega\times Y}A(G-H)\cdot(G-H)\,dxdy,
\]
even when $A$ is not symmetric. As $G-H=Qp-Q_mp_r$,
\[
 \lim_{\eps\to0}E_\eps
 \le\beta\|Q_mp_r-Qp\|_{L^2(\Omega\times Y)}^2 .
\]
Together with \eqref{eq:Eepsmr}, this proves \eqref{eq:energycorrector}. Finally, since $\|p_r\|_{L^2(\R^N)}\le\|p\|_{L^2(\R^N)}$,
\[
 \|Q_mp_r-Qp\|_{L^2(\Omega\times Y)}
 \le
 \|Q_m-Q\|_{L^2(Y)}\|p\|_{L^2(\R^N)}
 +\|Q\|_{L^2(Y)}\|p_r-p\|_{L^2(\R^N)},
\]
and \eqref{eq:doublelimitcorrector} follows.
\end{proof}

The corrected gradient must still be generated by an admissible function; this requires a cut-off near $\partial\Omega$. For $\delta>0$, choose $\eta_\delta\in C_c^\infty(\Omega)$ such that
\[
 0\le\eta_\delta\le1,
 \qquad
 \eta_\delta(x)=1\ \text{if }x\in\Omega,\ \dist(x,\partial\Omega)>2\delta,
\]
\[
 \supp\eta_\delta\subset\{x\in\Omega:\dist(x,\partial\Omega)>\delta\},
 \qquad
 |\nabla\eta_\delta|\le C\delta^{-1}.
\]
We also write $\Lambda_t:=\{x\in\R^N:\dist(x,\partial\Omega)<t\}$ for the two-sided boundary layer of width $t>0$. For fixed $m,r,\delta$, define
\[
 \mathcal C_\eps^{m,r,\delta}(x)
 :=\eps^s\eta_\delta(x)
 \sum_{j=1}^N\chi_{j,m}^s(x/\eps)(p_r(x))_j,\qquad x\in\R^N.
\]
Then $\mathcal C_\eps^{m,r,\delta}\in\Ho$, and Lemma~\ref{lem:amplitude}, applied to each trigonometric polynomial $\chi_{j,m}^s$ with amplitude $\eta_\delta(p_r)_j\in H^1(\R^N)$, gives
\begin{equation}\label{eq:admissiblerealization}
 \bigl\|
 \grads \mathcal C_\eps^{m,r,\delta}
 -\eta_\delta Q_m(x/\eps)p_r
 \bigr\|_{L^2(\R^N)}\le C_{m,r,\delta}\,\eps^s.
\end{equation}

\begin{theorem}\label{thm:generalcorrector}
Assume \eqref{eq:Aassumption}--\eqref{eq:ellipticity}. There exist functions
$m(\eps)\to\infty$, $r(\eps)\downarrow0$ and $\delta(\eps)\downarrow0$ such that
\[
 \mathcal C_\eps:=\mathcal C_\eps^{m(\eps),r(\eps),\delta(\eps)}\in\Ho
\]
satisfies
\begin{equation}\label{eq:Hscorrector}
 \|u_\eps-u_0-\mathcal C_\eps\|_{H^s(\R^N)}\to0
 \quad\text{and}\quad
 \bigl\|\grads\mathcal C_\eps-\eta_{\delta(\eps)}Q_{m(\eps)}(x/\eps)p_{r(\eps)}\bigr\|_{L^2(\R^N)}\to0 .
\end{equation}
\end{theorem}

\begin{proof}
Choose $\gamma_n\downarrow0$. First select $m_n\uparrow\infty$ and $r_n\downarrow0$ such that
\[
 \|Q_{m_n}-Q\|_{L^2(Y)}\|p\|_{L^2(\R^N)}+\|Q\|_{L^2(Y)}\|p_{r_n}-p\|_{L^2(\R^N)}<\gamma_n.
\] For these parameters choose $\delta_n\downarrow0$ such that $\|Q_{m_n}\|_{L^\infty(Y)}\|p_{r_n}\|_{L^2(\Omega\cap\Lambda_{3\delta_n})}<\gamma_n$. Only after fixing $m_n$, $r_n$ and $\delta_n$ do we choose $\eps_n\downarrow0$: by Proposition~\ref{prop:gradientcorrector} and \eqref{eq:admissiblerealization}, it can be chosen so that, for every $0<\eps\le\eps_n$, $\|\grads u_\eps-G_\eps^{m_n,r_n}\|_{L^2(\R^N)}\le C\gamma_n$ and $C_{m_n,r_n,\delta_n}\eps^s\le\gamma_n$. This order is essential because the realization constant $C_{m_n,r_n,\delta_n}$ need not remain bounded as $n\to\infty$. Choose the thresholds recursively with $0<\eps_n<\min\{\eps_{n-1}/2,1/n\}$ for $n\ge2$, and define $(m,r,\delta)(\eps)=(m_n,r_n,\delta_n)$ for $\eps_{n+1}<\eps\le\eps_n$. For $\eps>\eps_1$ set $\mathcal C_\eps=0$; this changes none of the asymptotic assertions. Since $\grads u_\eps-\grads u_0-\grads\mathcal C_\eps$ equals
\[
 \bigl(\grads u_\eps-G_\eps^{m,r}\bigr)+(\mathbf 1_\Omega-\eta_\delta)Q_m(x/\eps)p_r+\bigl(\eta_\delta Q_m(x/\eps)p_r-\grads\mathcal C_\eps\bigr),
\]
each term is at most $C\gamma_n$ in $L^2(\R^N)$. As $u_\eps-u_0-\mathcal C_\eps\in\Ho$, the fractional Poincar\'e inequality gives \eqref{eq:Hscorrector}.
\end{proof}

For the elliptic solutions, the global corrector estimate also controls the boundary layer, which the exterior compactness recalled in Subsection~\ref{sec:exterior} does not reach.

\begin{corollary}\label{cor:exteriorstrong}
Under \eqref{eq:Aassumption}--\eqref{eq:ellipticity}, $\grads u_\eps\to\grads u_0$ strongly in $L^2(\Omega^c;\R^N)$.
\end{corollary}

\begin{proof}
Since $\eta_{\delta(\eps)}Q_{m(\eps)}(x/\eps)p_{r(\eps)}$ vanishes on $\Omega^c$, the second convergence in \eqref{eq:Hscorrector} gives
\[
 \|\grads\mathcal C_\eps\|_{L^2(\Omega^c)}
 \le
 \bigl\|\grads\mathcal C_\eps-
 \eta_{\delta(\eps)}Q_{m(\eps)}(x/\eps)p_{r(\eps)}\bigr\|_{L^2(\R^N)}
 \longrightarrow0.
\]
The first convergence in \eqref{eq:Hscorrector} therefore implies
\[
 \|\grads u_\eps-\grads u_0\|_{L^2(\Omega^c)}
 \le\|\grads(u_\eps-u_0-\mathcal C_\eps)\|_{L^2(\R^N)}
 +\|\grads\mathcal C_\eps\|_{L^2(\Omega^c)}\longrightarrow0.
\]
\end{proof}

\begin{remark}\label{rem:correctorsHconv}
The nonlocal $H$-convergence theory \cite{BuchingerEtAl2026,CaponiCarbottiMaione2025} identifies the limit problem for general coefficients, not necessarily periodic, but it does not construct correctors. In the local theory, $H$-convergence comes with the correctors of Murat and Tartar \cite{Tartar2009}, which give strong convergence of the corrected gradients in $L^1_{\rm loc}$ but are not explicit; to our knowledge, no fractional analogue for general coefficients is available. Here the periodic structure yields explicit admissible correctors built from the fractional cell potentials, and Theorem~\ref{thm:generalcorrector} and Corollary~\ref{cor:exteriorstrong} give strong convergence in $H^s(\R^N)$, including near $\partial\Omega$ and in $\Omega^c$.
\end{remark}

\section{Homogenization with oscillatory flux forcing}\label{sec:oscforcing}

The elliptic two-scale theory also identifies effective source terms that are not determined by weak convergence of the right-hand sides alone. Fix $f\in(\Ho)'$, $g\in C_c^\infty(\Omega)$, and $F\in L^\infty_\#(Y;\R^N)$. Consider the variational problem
\begin{equation}\label{eq:oscforcingmicro}
 \int_{\R^N} A(x/\eps)\grads u_\eps\cdot\grads v\,dx
 =\langle f,v\rangle+\int_\Omega g(x)F(x/\eps)\cdot\grads v\,dx
 \qquad\forall v\in\Ho.
\end{equation}
Equivalently, $-\divs(A(x/\eps)\grads u_\eps)=f-\divs(gF(x/\eps))$ in $\Omega$, with the complementary-value condition. The last term is understood as a functional on $\Ho$, not as a pointwise divergence. Let $\vartheta_F^s\in\Hperzero$ solve
\begin{equation}\label{eq:oscforcingcell}
 \int_Y A(y)\gradys\vartheta_F^s\cdot\gradys\varphi\,dy
 =\int_Y F(y)\cdot\gradys\varphi\,dy
 \qquad\forall\varphi\in\Hperzero,
\end{equation}
and define the effective forcing vector
\[
 F^{\rm hom}:=\int_Y\bigl(F(y)-A(y)\gradys\vartheta_F^s(y)\bigr)\,dy.
\]
The cell problem is well posed by the Lax--Milgram theorem. By Theorem~\ref{thm:bridge}, $\vartheta_F^s=\mathcal J_s\vartheta_F$, where $\vartheta_F\in\dot H^1_\#(Y)$ solves the classical cell problem $\int_YA\nabla_y\vartheta_F\cdot\nabla_y\varphi\,dy=\int_YF\cdot\nabla_y\varphi\,dy$ for all $\varphi\in\dot H^1_\#(Y)$; hence $F^{\rm hom}$ coincides with the classical effective forcing and is independent of $s$. Proposition~\ref{prop:oscforcing} extends to the fractional setting the classical correction of oscillating divergence-form right-hand sides.

\begin{proposition}\label{prop:oscforcing}
Assume \eqref{eq:Aassumption}--\eqref{eq:ellipticity}. The solutions of \eqref{eq:oscforcingmicro} converge weakly in $\Ho$ and strongly in $L^2(\R^N)$ to the unique $u_0\in\Ho$ satisfying
\begin{equation}\label{eq:oscforcinglimit}
 \int_\Omega\Ahom\grads u_0\cdot\grads v\,dx
 +\int_{\Omega^c}\Abar\grads u_0\cdot\grads v\,dx
 =\langle f,v\rangle+\int_\Omega gF^{\rm hom}\cdot\grads v\,dx
 \qquad\forall v\in\Ho.
\end{equation}
Moreover, the interior fractional gradients satisfy
\begin{equation}\label{eq:oscforcing2scale}
 \grads u_\eps\twoweak
 [\Id+Q(y)]\grads u_0(x)+g(x)\gradys\vartheta_F^s(y)
 \quad\text{in }L^2(\Omega\times Y;\R^N).
\end{equation}
\end{proposition}

\begin{proof}
The bound $\|gF(x/\eps)\|_{L^2(\Omega)}\le\|g\|_2\|F\|_\infty$ gives a uniform bound for $u_\eps$ in $\Ho$, and the proof of Theorem~\ref{thm:homogenization} applies with the same realized tests $v_\eps=v_0+\eps^s v_1(x,x/\eps)$; only the forcing term is new. By Lemma~\ref{lem:periodicmultiplier}(i), applied to the constant family $ge_1$ and the matrix whose first column is $F$, $gF(x/\eps)\twoweak g(x)F(y)$ in $\Omega\times Y$. Writing $\grads v_\eps=\Psi^\eps+r_\eps$ with $\Psi:=\grads v_0+\gradys v_1\in L^2(\Omega;C_\#(Y;\R^N))$ and $r_\eps\to0$ in $L^2(\R^N)$, and testing with $\Psi$, we obtain
\[
 \int_\Omega gF(x/\eps)\cdot\grads v_\eps\,dx
 \longrightarrow
 \int_{\Omega\times Y}g(x)F(y)\cdot\Psi(x,y)\,dxdy.
\]
The two-scale identity \eqref{eq:twoscalevariational} thus holds with the additional term $\int_{\Omega\times Y}gF\cdot(\grads v_0+\gradys v_1)\,dxdy$ on the right-hand side. Taking $v_0=0$ and localizing in $x$, we obtain for almost every $x\in\Omega$ the periodic cell equation
\[
 \int_Y\bigl[A(y)(\grads u_0(x)+\gradys u_1(x,y))-g(x)F(y)\bigr]\cdot\gradys\varphi\,dy=0
 \qquad(\varphi\in\Hperzero).
\]
The uniqueness of \eqref{eq:cellfractional} and \eqref{eq:oscforcingcell} yields
\[
 u_1(x,y)=\sum_{j=1}^N(\grads u_0(x))_j\chi_j^s(y)+g(x)\vartheta_F^s(y),
\]
which proves \eqref{eq:oscforcing2scale}. Setting $v_1=0$ and averaging the two-scale identity in $y$ then yields \eqref{eq:oscforcinglimit}, initially for $v_0\in C_c^\infty(\Omega)$ and subsequently for all $v\in\Ho$ by density. The limit problem is coercive and hence uniquely solvable, so every subsequence has the same limit.
\end{proof}

\begin{remark}\label{rem:oscforcing}
Even if $\int_YF\,dy=0$, the vector $F^{\rm hom}$ need not vanish. In this case,
$gF(x/\eps)\rightharpoonup0$ in $L^2(\R^N;\R^N)$ (extended by
zero outside $\Omega$). By the boundedness of
$\divs:L^2(\R^N;\R^N)\to H^{-s}(\R^N)$, the induced
fractional-divergence forces converge only weakly in $(\Ho)'$ in general.
This weak convergence does not by itself identify the homogenized
right-hand side; the additional cell problem \eqref{eq:oscforcingcell}
records its correlation with the heterogeneous coefficient.
For example, when $N=1$, take $A=a$ scalar, with $a=2$ and $F=1$ on $(0,1/2)$,
and $a=1$ and $F=-1$ on $(1/2,1)$, extended periodically.
Then $\int_YF=0$, whereas the one-dimensional cell equation gives
\[
 F^{\rm hom}=\frac{\int_Y F(y)/a(y)\,dy}{\int_Y 1/a(y)\,dy}
 =\frac{1/4-1/2}{1/4+1/2}=-\frac13.
\]
\end{remark}

\section{Parabolic homogenization with dual-space forcing}\label{sec:parabolicsub}

We now let the data depend on time. Since the fast variable is purely spatial, the cell problems are unchanged. Under the assumptions of Theorem~\ref{thm:mainC}, $u_\eps\in L^2(0,T;H_0^s(\Omega))\cap H^1(0,T;(H_0^s(\Omega))')$ is the variational solution of
\begin{equation}\label{eq:paraboliceps}
 \left\langle\partial_tu_\eps,v\right\rangle
 +\mathfrak a_\eps(u_\eps,v)
 =\langle f_\eps,v\rangle
 \qquad\forall v\in H_0^s(\Omega),
\end{equation}
as an identity in $(H_0^s(\Omega))'$ for almost every $t\in(0,T)$, with $u_\eps(0)=u_\eps^{\rm in}$.

The space--time versions of compactness and of the multiplication lemmas read as follows. Two-scale convergence in $(0,T)\times\Omega\times Y$ is defined as in Subsection~\ref{sec:twoscale}, with test functions in $L^2((0,T)\times\Omega;C_\#(Y;\R^N))$ and the fast variable $x/\eps$ only.

\begin{lemma}\label{lem:paraboliccompactness}
Let $(w_\eps)$ be bounded in $L^2(0,T;H_0^s(\Omega))$ and suppose that $w_\eps\weakto w$ in that space. Up to a subsequence, there exists $w_1\in L^2((0,T)\times\Omega;\Hperzero)$, uniquely determined by the extracted two-scale gradient limit through its zero-mean normalization, such that
\[
 \grads w_\eps|_{(0,T)\times\Omega}\twoweak
 \grads w+\gradys w_1
 \quad\text{in }L^2((0,T)\times\Omega\times Y;\R^N).
\]
\end{lemma}

\begin{proof}
The proof of Theorem~\ref{thm:maincompactness} applies with $(t,x)\in(0,T)\times\Omega$ as slow variable. Let $G$ be a two-scale limit of $\grads w_\eps$. For tests $\varrho(t)\Phi(x,y)$ with $\varrho\in C_c^\infty(0,T)$ and $\Phi$ as in \eqref{eq:solconditions}, Lemma~\ref{lem:oscdiv} gives
\[
 \left|\int_0^T\!\!\int_\Omega\grads w_\eps\cdot\varrho(t)\Phi(x,x/\eps)\,dxdt\right|
 \le C\|\varrho\|_{L^2(0,T)}\sup_\eps\|w_\eps\|_{L^2((0,T)\times\Omega)}\,\eps^{1-s}.
\]
Countable dense families of tests in $(t,x)$ and of smooth periodic solenoidal fields then yield
\[
 G(t,x,\cdot)-\grads w(t,x)\perp\mathcal Z_s
 \qquad\text{for a.e.\ }(t,x)\in(0,T)\times\Omega.
\]
The Fourier inverse \eqref{eq:potentialfourier} defines $w_1$ by a bounded linear map of $G-\grads w$, which gives Bochner measurability and
\[
 \|w_1\|_{L^2((0,T)\times\Omega;\Hperzero)}=\|G-\grads w\|_{L^2((0,T)\times\Omega\times Y)}.
\]
\end{proof}

\begin{lemma}\label{lem:parabolicmultiplier}
Let $B\in L^\infty_\#(Y;\R^{N\times N})$ and $B_\eps(x)=B(x/\eps)$.
\begin{enumerate}[label=(\roman*)]
\item If $(z_\eps)$ is bounded in $L^2((0,T)\times\Omega;\R^N)$ and $z_\eps\twoweak z$, then
\[
 B_\eps z_\eps\twoweak B(y)z(t,x,y).
\]
\item If $w_\eps\to w$ strongly in $L^2((0,T)\times\Omega)$ and $(w_\eps)$ is bounded in $L^2(0,T;H_0^s(\Omega))$, then for every $h\in L^2((0,T)\times\Omega^c;\R^N)$,
\begin{equation}\label{eq:parabolicexteriorpairing}
 \int_0^T\!\!\int_{\Omega^c}B_\eps\grads w_\eps\cdot h\,dxdt
 \longrightarrow
 \int_0^T\!\!\int_{\Omega^c}\overline B\,\grads w\cdot h\,dxdt.
\end{equation}
\end{enumerate}
\end{lemma}

\begin{proof}
(i) The proof of Lemma~\ref{lem:periodicmultiplier}(i) applies verbatim, with $(t,x)\in(0,T)\times\Omega$ as slow variable and the oscillation acting only on $x$.

(ii) For a.e.\ $t$, $w_\eps(t)-w(t)$ vanishes outside $\overline\Omega$. For every open $\Omega_1\Supset\Omega$, Lemma~\ref{lem:farfield} therefore gives
\[
 \grads w_\eps\to\grads w\quad\text{in }L^2((0,T)\times\Omega_1^c;\R^N),
\]
while Lemma~\ref{lem:periodicmultiplier}(ii) at a.e.\ $t$ and dominated convergence give
\[
 B_\eps^Th\weakto\overline B^{\,T}h\quad\text{in }L^2((0,T)\times\Omega_1^c;\R^N).
\]
The rest of the proof of Lemma~\ref{lem:exteriorpairing} applies with these two facts in place of the compact embedding and of Lemma~\ref{lem:periodicmultiplier}(ii); the layer $\Omega_1\setminus\Omega$ is chosen, using $|\partial\Omega|=0$, so that $\|h\|_{L^2((0,T)\times(\Omega_1\setminus\Omega))}$ is small.
\end{proof}

With these tools we can pass to the limit in \eqref{eq:paraboliceps}.

\begin{theorem}\label{thm:parabolic}
Under the assumptions of Theorem~\ref{thm:mainC}, the solutions of \eqref{eq:paraboliceps} satisfy
\begin{align}
 u_\eps&\weakto u_0
 &&\text{in }L^2(0,T;H_0^s(\Omega)),\nonumber\\
 u_\eps&\to u_0
 &&\text{in }L^2((0,T)\times\Omega),\label{eq:parstrong}
\end{align}
and
\[
 \grads u_\eps
 \twoweak
 \grads u_0+\gradys u_1
 \quad\text{in }L^2((0,T)\times\Omega\times Y;\R^N),
 \qquad
 u_1(t,x,y)
 =\sum_{j=1}^N(\grads u_0(t,x))_j\chi_j^s(y).
\]
Moreover,
\[
 u_0\in L^2(0,T;H_0^s(\Omega))\cap H^1\bigl(0,T;(H_0^s(\Omega))'\bigr)
 \subset C([0,T];L^2(\Omega)),
\]
and $u_0$ is the unique variational solution of
\begin{equation}\label{eq:parhom}
 \left\langle\partial_tu_0,v\right\rangle
 +\mathfrak a_0(u_0,v)
 =\langle f,v\rangle
 \quad\forall v\in H_0^s(\Omega),\ \text{for a.e.\ }t\in(0,T),
 \qquad
 u_0(0)=u^{\rm in}.
\end{equation}
\end{theorem}

\begin{proof}
Set $V=H_0^s(\Omega)$ and $V'=(H_0^s(\Omega))'$. Since $\mathfrak a_\eps$ is bounded and uniformly $V$-coercive, the Gelfand-triple evolution theorem gives a unique solution in $L^2(0,T;V)\cap H^1(0,T;V')\subset C([0,T];L^2(\Omega))$. Testing with $u_\eps$, and using the fractional Poincar\'e and Young inequalities and then the equation, we obtain, uniformly in $\eps$,
\[
 \|u_\eps\|_{L^\infty(0,T;L^2(\Omega))}+\|u_\eps\|_{L^2(0,T;V)}+\|\partial_tu_\eps\|_{L^2(0,T;V')}
 \le C\bigl(\|u_\eps^{\rm in}\|_{L^2(\Omega)}+\|f_\eps\|_{L^2(0,T;V')}\bigr).
\]
By fractional Rellich compactness, $V\Subset L^2(\Omega)\hookrightarrow V'$, and \cite[Corollary~4]{Simon1987}, applied to this triple, yields \eqref{eq:parstrong} after extraction; the time derivatives converge weakly in $L^2(0,T;V')$ to $\partial_tu_0$.

Lemma~\ref{lem:paraboliccompactness} provides a space--time two-scale limit
\[
 \grads u_\eps\twoweak G:=\grads u_0+\gradys u_1,
 \qquad
 u_1\in L^2((0,T)\times\Omega;\Hperzero),
\]
and Lemma~\ref{lem:parabolicmultiplier}(i) gives $A_\eps\grads u_\eps\twoweak A(y)G$. By Lemma~\ref{lem:parabolicmultiplier}(ii) and \eqref{eq:parstrong}, the exterior pairing \eqref{eq:parabolicexteriorpairing} holds with $B=A$.

To pass to the limit we use realized oscillatory test functions. Let $\psi_0\in C_c^\infty([0,T)\times\Omega)$, let $\psi_1\in C_c^\infty((0,T)\times\Omega;C_\#^\infty(Y))$ have zero $Y$-mean, and set
\[
 \psi_1(t,x,y)=\sum_{k\ne0}a_k(t,x)e^{iq_k\cdot y},
 \qquad
 \psi_\eps(t,x)=\psi_0(t,x)+\eps^s\psi_1(t,x,x/\eps).
\]
Applying the realization estimate \eqref{eq:strongrealization} at almost every $t$ and then taking the $L^2(0,T)$ norm gives
\[
 \grads\psi_\eps=\Psi^\eps+r_\eps,
 \qquad
 \Psi:=\grads\psi_0+\gradys\psi_1,
\]
with
\[
 \|r_\eps\|_{L^2((0,T)\times\R^N)}
 \le C_s\eps^s
 \Bigl\|\sum_{k\ne0}\|a_k(t,\cdot)\|_{\dot H^s(\R^N)}\Bigr\|_{L^2(0,T)}
 \le C_{s,\psi_1}\eps^s.
\]
The Fourier-amplitude norm is finite because $\psi_1$ is smooth and compactly supported in $(t,x)$, and in the interior $\Psi$ belongs to $L^2((0,T)\times\Omega;C_\#(Y;\R^N))$. Outside $\Omega$, the realized microscopic term vanishes, while $\grads\psi_0$ may be nonzero because the fractional gradient is nonlocal. Moreover,
\[
 \partial_t\psi_\eps\to\partial_t\psi_0\ \text{in }L^2((0,T)\times\Omega),
 \qquad
 \psi_\eps(0)=\psi_0(0),
 \qquad
 \psi_\eps\weakto\psi_0\ \text{in }L^2(0,T;H_0^s(\Omega)),
\]
the last convergence because the family is bounded in $L^2(0,T;H_0^s(\Omega))$ and converges in $L^2((0,T)\times\Omega)$.

The time-integrated weak formulation reads
\[
 -\int_0^T(u_\eps,\partial_t\psi_\eps)\,dt+\int_0^T\mathfrak a_\eps(u_\eps,\psi_\eps)\,dt
 =\int_0^T\langle f_\eps,\psi_\eps\rangle\,dt+(u_\eps^{\rm in},\psi_0(0)).
\]
The first term converges by \eqref{eq:parstrong}, and the right-hand side by the strong convergence of $f_\eps$ and $u_\eps^{\rm in}$ and the weak convergence of $\psi_\eps$. In $\mathfrak a_\eps$, the space--time two-scale convergence of $A_\eps\grads u_\eps$ tested with $\Psi$, together with $\|r_\eps\|_{L^2((0,T)\times\R^N)}\to0$, and the exterior pairing give
\[
 \int_0^T\mathfrak a_\eps(u_\eps,\psi_\eps)\,dt
 \longrightarrow
 \int_0^T\!\!\int_{\Omega\times Y}A(y)G\cdot\Psi\,dxdydt
 +\int_0^T\!\!\int_{\Omega^c}\Abar\grads u_0\cdot\grads\psi_0\,dxdt.
\]
Hence, for all such $(\psi_0,\psi_1)$,
\begin{equation}\label{eq:partwoscalevar}
 \begin{aligned}
 &-\int_0^T(u_0,\partial_t\psi_0)
 +\int_0^T\!\!\int_{\Omega\times Y}A(y)G\cdot(\grads\psi_0+\gradys\psi_1)
 +\int_0^T\!\!\int_{\Omega^c}\Abar\grads u_0\cdot\grads\psi_0\\
 &\qquad=\int_0^T\langle f,\psi_0\rangle+(u^{\rm in},\psi_0(0)).
 \end{aligned}
\end{equation}

To identify the limit, take $\psi_0=0$ and $\psi_1=\phi(t,x)\varphi(y)$, with $\phi\in C_c^\infty((0,T)\times\Omega)$ and $\varphi\in C_\#^\infty(Y)$ of zero mean. This gives, for a.e.\ $(t,x)\in(0,T)\times\Omega$, the cell problem \eqref{eq:cellfractional} with $\xi=\grads u_0(t,x)$, hence
\[
 u_1(t,x,\cdot)=\sum_{j=1}^N(\grads u_0(t,x))_j\chi_j^s,
 \qquad
 \int_YA(y)G\,dy=\Ahom\grads u_0.
\]
Taking then $\psi_1=0$ and $\psi_0(t,x)=\varrho(t)v(x)$, with $\varrho\in C_c^\infty((0,T))$ and $v\in C_c^\infty(\Omega)$, gives \eqref{eq:parhom} in the sense of distributions on $(0,T)$, hence a.e.\ in $t$ by density in $v$.

To identify the initial value, take now $\psi_1=0$ and $\psi_0(t,x)=\varrho(t)v(x)$ with $v\in C_c^\infty(\Omega)$ and $\varrho\in C_c^\infty([0,T))$. The limit equation, together with the Gelfand-triple integration-by-parts formula, gives
\[
 -\int_0^T(u_0(t),v)_{L^2(\Omega)}\varrho'(t)\,dt
 =\int_0^T\langle\partial_tu_0(t),v\rangle\varrho(t)\,dt
   +(u_0(0),v)_{L^2(\Omega)}\varrho(0).
\]
On substituting $\langle\partial_tu_0,v\rangle=\langle f,v\rangle-\mathfrak a_0(u_0,v)$ and comparing with \eqref{eq:partwoscalevar}, all integral terms cancel and leave
\[
 (u_0(0)-u^{\rm in},v)_{L^2(\Omega)}\,\varrho(0)=0.
\]
Choosing $\varrho(0)=1$ and using the density of $C_c^\infty(\Omega)$ in $L^2(\Omega)$ yields $u_0(0)=u^{\rm in}$ in $L^2(\Omega)$. The limit problem \eqref{eq:parhom} has a unique solution, so the whole family converges.
\end{proof}

\section*{Declarations}

\noindent\textbf{Data availability.} No data were used for the research described in this article.

\smallskip
\noindent\textbf{Funding.} This research did not receive any specific grant from funding agencies in the public, commercial, or not-for-profit sectors.

\smallskip
\noindent\textbf{Conflict of interest.} The author declares that he has no conflict of interest.

\end{document}